\documentclass[ejs]{imsart}

\RequirePackage[OT1]{fontenc}
\RequirePackage{amsthm,amsmath, amssymb, bbm, epsfig, booktabs, url, color}
\RequirePackage[numbers]{natbib}
\RequirePackage[colorlinks,citecolor=blue,urlcolor=blue]{hyperref}

\startlocaldefs
\numberwithin{equation}{section}
\theoremstyle{plain}
\newtheorem{thm}{Theorem}[section]
\newtheorem{cor}{Corollary}[section]
\newtheorem{lem}{Lemma}[section]
\newtheorem{prop}{Proposition}[section]
\newtheorem{rem}{Remark}[section]
\newtheorem{ex}{Example}[section]
\endlocaldefs

\newcommand{\beq}{\begin{equation}}
\newcommand{\eeq}{\end{equation}}
\newcommand{\beqa}{\begin{eqnarray}}
\newcommand{\eeqa}{\end{eqnarray}}
\newcommand{\beqas}{\begin{eqnarray*}}
\newcommand{\eeqas}{\end{eqnarray*}}

\newcommand{\ceq}{\!\!\! & = & \!\!\!}
\newcommand{\cleq}{\!\!\! & \le & \!\!\!}

\newcommand{\cF}{\mathcal{F}}
\newcommand{\cG}{\mathcal{G}}

\newcommand{\cN}{\mathcal{N}}

\newcommand{\real}{{\ensuremath{\mathbb{R}}}}

\newcommand{\el}{{\epsilon}}

\newcommand{\sumnl}{{\sum\nolimits}}
\newcommand{\supnl}{{\sup\nolimits}}

\newcommand{\maxnl}{{\max\nolimits}}

\newcommand{\limnl}{{\lim\nolimits}}
\newcommand{\one}{\mathbbm{1}}

\newcommand{\done}{\hfill $\Box$}

\definecolor{green}{rgb}{0.1,0.55,0.1}
\definecolor{blau}{rgb}{0,0.4,1}

\begin{document}

\begin{frontmatter}
	
	\title{Forecast dominance testing via sign randomization}

\begin{aug}
\author{\fnms{Werner} \snm{Ehm}\ead[label=e1]{wernehm@web.de}}

\address{Heidelberg Institute for Theoretical Studies\\ Computational Statistics Group \\ Schloss-Wolfsbrunnenweg 35\\ 69118 Heidelberg, Germany\\ 
\printead{e1}}

\author{\fnms{Fabian } \snm{Kr\"uger}
\ead[label=e2]{fabian.krueger@awi.uni-heidelberg.de}}

\address{Heidelberg University\\ Alfred-Weber-Institute for Economics\\ Bergheimer Strasse 58 \\ 69115 Heidelberg, Germany\\
\printead{e2}}

\runauthor{W. Ehm \and F. Kr\"uger}
\runtitle{Randomization Test}

\affiliation{Some University and Another University}

\end{aug}

\begin{abstract}
We propose randomization tests of whether forecast 1 outperforms forecast 2 across a class of scoring functions. This hypothesis is of applied interest: While the prediction context often prescribes a certain class of scoring functions, it is typically hard to motivate a specific choice on statistical or substantive grounds. We investigate the asymptotic behavior of the test statistics under mild conditions, avoiding the need to assume particular dynamic properties of forecasts and realizations. 
The properties of the one-sided tests depend on a corresponding version of Anderson's inequality, which we state as a conjecture of independent interest. 
Numerical experiments and a data example indicate that the tests have good size and power properties in practically relevant situations. 
\end{abstract}

\begin{keyword}[class=MSC]
\kwd{62G09}
\kwd{62G10}
\kwd{62M20}
\end{keyword}

\begin{keyword}
\kwd{comparative forecast evaluation}
\kwd{hypothesis testing}
\kwd{randomization}
\end{keyword}
\tableofcontents
\end{frontmatter}

\section{Introduction}
\label{sec:intro}

Forecasts of future events and quantities are essential across disciplines. At the same time, forecasts notoriously are imprecise and prone to bias, calling for methods to assess and compare the performance of imperfect predictions both theoretically and on the basis of empirical data. In the context of point forecasts, which we consider in this paper, the appropriate evaluation tool is that of a {\em consistent scoring function} \cite{G}. 
A scoring function $S\equiv S(x,y)$ assigns to each forecast $x$ and realization $y$ a real-valued score such that smaller scores correspond to better forecasts. Specifically, let $\phi$ be a real-valued functional defined on a class $\cG$ of possible distributions $G$ of $y$, such as the mean or a quantile of $G$. The scoring function is {\em consistent (for $\phi$ relative to $\cG$)} if $S(\phi(G),G) \le S(x,G)$ for every $G \in \cG$ and forecast $x$; here $S(x,G)= E_{y \sim G}~S(x, y)$. With a consistent scoring function the forecaster can do no better than predict the true functional value, which rewards honest reporting. 

For a given functional $\phi$, supposed fixed in the following, there generally exists a whole class of consistent scoring functions, or ``scores.'' Characterizations of the respective score classes for various functionals may be found in e.g.~\cite{BGW,FZ,G}. For example, all scores of the form $S(x, y) = \varphi(y) - \varphi(x) - \varphi'(x)~(y-x)$, $\varphi$ a convex function with subgradient $\varphi'$, are consistent for the mean functional \cite{BGW,S}. 
In applied contexts, consistent alternatives to the special case $\varphi(x) = x^2$ of the squared error score were discussed in \cite{GRaf}, \cite{MS}, \cite{SAM} and others for the binary case $y \in \{0, 1\}$, and in \cite{Pt0} for positive predictands $y \in \mathbb{R}_+$. 

The availability of an entire family of scoring functions that are theoretically legitimate comes with the drawback that two scores may produce different forecast rankings even if both are consistent for the same target functional \cite{EGJK,EGK,Pt}. This lack of robustness is unsatisfactory as there are often no strong arguments for choosing a particular score. It is therefore natural to ask whether forecast rankings are stable across a family of scores. E.g., given two forecasts $x_1, x_2 ,\ $ does $x_1$ dominate $x_2$ in the sense that it is superior with respect to {\em all} consistent scoring functions? In the case of a quantile or expectile functional\footnote{Expectiles are an asymmetric generalization of the mean; they were introduced by \cite{NP} and have recently received attention in financial risk management. We provide a formal definition of expectiles in Section \ref{sec:asy}.} the complexity of the problem can be much reduced by means of a Choquet representation: here every consistent score $S$ can be represented as a mixture of ``elementary'' (in fact: extremal) scores $S_\theta,\, \theta \in \real$. That is, for every $S$ there is a nonnegative Borel measure $M$ on $\real$ such that $S(x,y) = \int S_\theta(x,y)\,  dM(\theta)$ \cite{EGJK}. The form of the elementary scores depends on the specific functional being studied; for example, in case of the mean functional 
\beq\label{escm}
S_\theta(x, y) = \begin{cases} |y-\theta|\ &\ \text{if \ min}\{x, y\} \equiv x\wedge y \le \theta < x\vee y \equiv \text{max}\{x, y\}, \\ 0 &\  \text{else}. \end{cases}
\eeq
The Choquet representation makes it possible to reduce dominance with respect to all consistent scoring functions to dominance with respect to the linearly indexed family of elementary scores $\{S_\theta:\, \theta \in \real\}$, a substantial simplification. Nevertheless, testing related hypotheses remains a challenging task.

Our focus will be on a very simple test of forecast dominance that goes without any of the common assumptions. In return, the notion of forecast dominance acquires a central role. Initially, it may be defined as follows \cite{EGJK}. Given a family $\mathcal{S} = \{S_\theta:\, \theta \in \Theta\}$ of (consistent) scoring funcions $S_\theta$ we say that {\em forecast $x_1$ dominates forecast $x_2$ at the distribution $G \in \cG$ if $S_\theta(x_1,G) \le S_\theta(x_2,G)$ for every} $\theta \in \Theta$. In such a one-step scenario, the set of all $G \in \cG$ satisfying this condition could constitute the hypothesis `$x_1$ dominates $x_2$' (with respect to $\mathcal{S}$). The question how to formulate suitable dominance hypotheses becomes substantially more involved when, as usually, forecasts $x_{k1},\, x_{k2}$ are produced step by step and the realizations $y_k$ become known before the next forecast instance. Current work on forecast evaluation and comparison emphasizes the joint dynamic behavior of forecasts and realizations, by using martingale methods \cite{LGS,SD}, the concept of prediction spaces \cite{GR,SZ}, or comparisons of conditional predictive ability \cite{GW}. Bootstrap tests of dominance hypotheses are presented in, e.g., \cite{JCS,LMW,LSW,Yen}, based on different dominance concepts and asymptotic frameworks. An account of the related literature addressing the relations with, and differences to the present approach is given in the discussion section \ref{sec:disc}.

Usually, mathematical analyses proceed from statistical models for the data and the formulation of hypotheses to related tests and their properties. Here we follow a reverse path. We make no assumptions about possible data generating mechanisms; instead we focus on a simple-to-implement test procedure and {\em ask for hypotheses for which this procedure represents a valid test} (asymptotically, at a given level). 
We take this route because quite often very little is known about the stochastic nature of the data. In fact, typical forecasting problems have to cope with complex statistical dependencies, structural change, and limited domain knowledge. Thus presumably, most of the usual assumptions do not apply, with largely unknown consequences, and are hard or impossible to check. We therefore have recourse to the classical Fisherian technique of external randomization, which is completely under one's control, and treat everything conditionally on the data $(x_{k1},\, x_{k2},\, y_k),\,k = 1,\ldots,n$. 

The use of external randomization to compare forecast performance dates back at least to \cite{DM}. Here we compare forecast performance across families of scores, rather than with respect to a single scoring function. 
Concretely, our goal is to elaborate on the sign randomization procedure tentatively proposed in \cite[end of Section 3]{EGJK} for testing forecast dominance. The idea is to reject the hypothesis `forecast 1 dominates forecast 2' if, e.g., $\sup_{\theta\in \Theta}\,  D_n(\theta)$ exceeds some critical value $c_n$ where 
$$
D_n(\theta) = n^{-1/2}\, \sumnl_{k \le n} d_k(\theta),\quad d_k(\theta) =  S_\theta(x_{k1},y_k) - S_\theta(x_{k2},y_k) 
$$
is the sum of the single score differences $d_k(\theta)$, properly scaled for the sake of asymptotics. Unfortunately, determination of $c_n$ generally is very diffcult even asymptotically; it appears impossible without making assumptions about the stochastic structure of the data. Our suggestion in \cite{EGJK} was to determine $c_n$ such that $Pr^* [\sup_{\theta\in \Theta}\,  D_n^*(\theta) > c_n] \approx \alpha$, the test level, where 
$$
D_n^*(\theta) = n^{-1/2}\, \sumnl_{k \le n} d_k(\theta) \sigma_k
$$
and $Pr^*$ exclusively refers to the i.i.d.~(``Rademacher'') random variables $\sigma_k$ assuming the values $\pm 1$ with probability 1/2 each. This clearly raises questions. 

First, how can the randomization distribution be connected to the distribution of the test statistic, particularly when no model assumptions are being made? Secondly, what precisely is to be understood under the hypothesis `forecast 1 dominates forecast 2'? As explained in Section \ref{sec:heu}, there is in fact a close connection between the two problems that helps to get around both -- up to one missing link: The approximate validity of our one-sided tests depends on an unproven variant of the celebrated Anderson's inequality \cite{A}. While for symmetric hypotheses postulating `no difference in predictive performance' the classical Anderson's inequality provides the necessary link, the asymmetry of dominance hypotheses requires a one-sided version of the inequality which we state as a conjecture that appears of independent interest. 

Obviously, dispensing with model asssumptions cannot mean doing without any assumptions. However, as detailed in Section \ref{sec:condwc}, asssumptions distantly related to stationarity and (in-)dependence properties of forecasts and observations will enter in a very indirect manner only, via basic asymptotic stability and ``moderate local clustering'' conditions, respectively, which are fulfilled under virtually any of the standard statistical models; cf.~Section \ref{sec:addmat}. Building on this novel asymptotic framework we present, in Section \ref{sec:asy}, weak convergence results governing the asymptotics of our test statistics in the important special cases of quantile and expectile forecasts.
For the overall organization of the paper see the table of contents. R \cite{R} program code to implement the randomization test is available at \url{https://github.com/FK83/fdtest}.

\section{Testing for forecast dominance -- Initial considerations}
\label{sec:ic}

\subsection{Formal setup}
\label{ss:setup}

Let $(x_{k1},x_{k2},y_k),\, k = 1,\ldots,n$ be a sequence of $n$ triplets where $x_{k1},x_{k2}$ are two point forecasts each for the subsequent observation $y_k$. The triplets are considered as random variables on a common probability space $(\Omega,\cF,Q)$ endowed with a filtration $\{\cF_k, k=0,\ldots,n\}$ such that $(x_{k1},x_{k2},y_k)$ is $\cF_k$-measurable for every $k$, and $\cF_0$ is trivial. Given a family $\mathcal{S} = \{S_\theta,\, \theta \in \Theta\}$ of scoring functions $S_\theta$, we compare the two forecasts via the suitably normalized average difference of their $S_\theta$ scores, i.e., we are interested in the stochastic process 
\beq\label{Dndef}
\theta \mapsto D_n(\theta) = n^{-1/2}\, \sumnl_{k\le n}\, d_k(\theta), 
\quad  d_k(\theta) = S_\theta( x_{k1},y_k) - S_\theta( x_{k2},y_k).
\eeq
Initially, the family $\mathcal{S}$ may be arbitrary; from Section \ref{sec:asy} on it will be specialized to the elementary scores for quantile and expectile forecasts.

\subsection{Notions of forecast dominance}
\label{sec:nfd}

Which notion of forecast dominance fits with the randomizarion test in our focus? 
One possibility to introduce forecast dominance in the present framework is to declare forecast 1 as weakly dominating forecast 2 at $Q$ (with respect to $\mathcal{S}$) if $E_Q\, D_n(\theta) \le 0$ for every $\theta\in \Theta$. The same condition furnishes a natural one-sided hypothesis in a testing context: 
\beq\label{Hweak}
H_-^w: \quad \supnl_\theta\ E_Q\, D_n(\theta) = \supnl_\theta\  n^{-1/2}\, \sumnl_{k\le n}\, E_Q\, d_k(\theta) \le 0.
\eeq
In fact, $H_-^w$ stands for all probability measures $Q$ under which $\sup\nolimits_\theta E_Q\, D_n(\theta) \le 0$. Note that the hypothesis does not depend on the scaling constant, which could also be taken as $n^{-1}$, or the constant 1. The present choice, $n^{-1/2}$, is simply for convenience in regard to the large sample asymptotics to be discussed later on.

The hypothesis $H_-^w$ involves unconditional expectations referring to both the observations $y_k$ and the forecasts $x_{k\ell}$. This form of $H_-^w$ is at odds with the dominance concept of Section \ref{sec:intro} which is sequential in nature and makes no assumption about the dynamics of the forecasts, hence is more flexible in this sense. Better accordance with this initial concept is achieved on replacing the unconditional expectations $E_Q\, d_k(\theta)$ in (\ref{Hweak}) by conditional expectations given the past. This leads upon the following definition of forecast dominance: we say that {\em forecast 1 dominates forecast 2 at $Q$} if
\beq\label{Mneg}
M_{n,Q}(\theta) = n^{-1/2}\,\sumnl_{k\le n}\, E_Q\, [\, d_k(\theta)\, |\, \cF_{k-1}] \le 0 \quad (\mbox{$Q$-a.s.,} \ \theta\in \Theta)
\eeq
(a.s.~is short for almost surely).
The corresponding hypothesis, $H_-$, then comprises all probabilities $Q$ for which (\ref{Mneg}) holds,
\beq\label{H_def}
\!\!\!\!\!\!\!\!\!\!\!\! H_-\, :\quad\, M_{n,Q}(\theta)\le 0 \quad (\mbox{$Q$-a.s.,} \ \theta\in \Theta).
\eeq
Hypothesis $H_-$ is defined by conditions on a random process rather than on a parameter, which is unusual. Nevertheless, it makes for a useful notion of forecast dominance, and it will prove to be the appropriate hypothesis for the randomization test. We also will consider sub-hypotheses of $H_-$ including, specifically, the hypothesis
\beq\label{H_sdef}
H_-^s\, :\quad\, E_Q\, [\, d_k(\theta)\, |\, \cF_{k-1}] \le 0 \quad (\mbox{$Q$-a.s.,} \ \theta\in \Theta,\ k=1,\ldots,n).\footnote{If the functions $\theta \mapsto d_k(\theta)$ are sufficiently regular, the conditions (\ref{H_def}), (\ref{H_sdef}) hold $Q$-a.s.~{\em for all $\theta \in \Theta$ simultaneouls}y. We do not further dwell on this technicality.}
\eeq
The interpretation of $H_-^s$ is straightforward: It says that forecast 1 is at least as good as forecast 2 at each time step. 
In \cite{GW}, this dominance concept is referred to as a comparison of {\em conditional predictive ability.}  In this parlance, the hypotheses $H_-$ and $H_-^s$ express superiority of forecast 1 over forecast 2 in terms of, respectively, average and step-by-step conditional predictive ability.

\begin{ex}\label{xmpl}
{\em For illustration we consider the case where the forecasters know one part each of the verifying observation. Specifically, let $y_k = \eta_{k1} + \eta_{k2}$ where $\eta_{k\ell},\, k \ge 1,\, \ell = 1,2$ are two independent autoregressive processes of the form $\eta_{k\ell} = a\, \eta_{k-1,\ell} + \el_{k\ell}$ with the same parameter $a$ and independent innovations $\el_{k\ell} \sim \cN(0,\tau_\ell^2)$. Suppose that at any instance $k$, forecaster 1 has access to $\eta_{k1}$ and the preceding value $\eta_{k-1,2}$ 
of the second process.\footnote{This setup is similar to the simulation example in \cite[Section 4.1]{GR}, except that our variant includes time series dynamics.} If  $a$ is known, a natural choice for the prediction of $y_k$ is $x_{k1} = \eta_{k1} + a\,\eta_{k-1,2}$. By definition,
$$
x_{k1} = a\,\eta_{k-1,1} + \el_{k1} + a\, \eta_{k-1,2} = a\, y_{k-1} + \el_{k1},
$$
and if forecaster 2's prediction similarly is $x_{k2} = \eta_{k2} + a\, \eta_{k-1,1}$, then $x_{k2} = a\, y_{k-1} + \el_{k2}$, and of course, $y_k = a\, y_{k-1} + \el_{k1} + \el_{k2}$. Taking at first squared error as the scoring function, the $k$-th score difference becomes
$$
d_k = (y_k- x_{k1})^2 - (y_k- x_{k2})^2 = \el_{k2}^2 - \el_{k1}^2.
$$
Thus if $\tau_1 \ge \tau_2$, say, and the innovations $\el_{k\ell}$ are independent of $\cF_{k-1}$ -- as in the common case where $\cF_k$ is the $\sigma$-algebra generated by all triplets $(x_{j1}, x_{j2}, y_j),\, j \le k$ --, then $E \, [\, d_k\, |\, \cF_{k-1}] \le 0$, consistent with the intuition that the forecaster having access to the more variable component should be better off. 
The case $\tau_1 = \tau_2$ is an instance of a situation where, a priori, none of the two forecasters is believed to outperform the other. Here $E \, [\, d_k\, |\, \cF_{k-1}] = 0$, which holds in fact for {\em every} scoring function $S$ within the present model: by symmetry the joint conditional distributions of $(x_{kj},y_k)=(a\, y_{k-1}+\el_{kj}, a\, y_{k-1}+\el_{k1}+\el_{k2})$ given $\cF_{k-1}$ are identical $(j=1,2)$, whence the conditional expectation of $d_k = S(x_{k1},y_k) - S(x_{k2},y_k)$ vanishes. 
In particular, $E\, [d_k(\theta) \mid \cF_{k-1}] = 0$ for every $\theta\in \real$, where $d_k(\theta)$ denotes the score difference with respect to the elementary scoring functions $S_\theta$ for the mean value; see (\ref{escm}). In this case we even have the following. 
}
\end{ex}

\begin{prop}\label{prop1}
If $\tau_1 > \tau_2$ then $E\, [d_k(\theta) \mid \cF_{k-1}] \le 0$  for every $\theta\in \real$ and $k\le n$, i.e., $H_-^s$ (and a fortiori $H_-$) holds.
\end{prop}

The example indicates that $H_-$ and $H_-^s$ represent meaningful conditions characterizing different forms of forecast dominance. Further examples may be found in the recent systematic study \cite{KrZi} of dominance conditions based on the concept of convex order.

\subsection{A fictitious test}

For testing the one-sided hypothesis $H_-$ (or $H_-^s$) it is natural to work with test functionals $T = T(D_n)$ that are {\em monotone} in the sense that $T(f) \le T(g)$ whenever $f \le g$ pointwise on $\Theta$, an interval in $\real$ from now on. Main examples are the supremum functional and integrals of the positive part of $D_n$, e.g.
\beq\label{Tdef}
T_\infty(f) = \supnl_{\theta\in \Theta}\, f(\theta), \quad T_p(f) = \int_\Theta f(\theta)_+^p d\theta\quad (p \ge 1,\ a_+ = a\vee 0).
\eeq
These functionals are convex in $f$, and accordingly, {\em the generic test functional $T$ is supposed to be monotone and convex} in the real functions $f$ on $\Theta$.

Importantly for the following, in tests of $H_-$ or $H_-^s$ based on such a functional $T$ it suffices to control the error of the first kind at the ``boundary'' of the hypothesis, where either $M_{n,Q}(\theta) \equiv 0$ or $E_Q\, [d_k(\theta) \mid \cF_{k-1}] \equiv 0$. Indeed, put 
$$
\widetilde d_{k,Q}(\theta) = d_k(\theta) - E_Q\, [\, d_k(\theta)\, |\, \cF_{k-1}],
$$
and let
$$
\widetilde D_{n,Q}(\theta) = n^{-1/2}\,  \sumnl_{k \le n}\, \widetilde d_{k,Q}(\theta) = D_n(\theta) - M_{n,Q}(\theta)
$$
denote the conditionally centered version of $D_n$. Then by monotonicity
\beqa
\!\!\!\!\!\!\!\!\!\! Pr_Q\, [T(D_n) > c] \ceq Pr_Q\, [T(\widetilde D_{n,Q} + M_{n,Q}) > c] \le 
Pr_Q\, [T(\widetilde D_{n,Q}) > c] \label{alfacntrl}
\eeqa
for every $Q \in H_-$, as claimed. 
Thus, {\em if critical values $\widetilde c_n$ could be obtained such that $\sup_{Q \in {H_-}}\! Pr_Q\, [T(\widetilde D_{n,Q}) > \widetilde c_n] \approx \alpha$,\footnote{The notation `$a_n \approx b_n$' is short for $a_n-b_n \to 0$ as $n\to \infty$.}\, the rule `reject $H_-$ if $T(D_n) > \widetilde c_n$' would give us an approximate level-$\alpha$ test of $H_-$.} 
For compactness of notation we henceforth drop the index $Q$, taking the dependence on the underlying probability measure as self-evident.

An initial step toward the determination of critical values is the following proposition, for which we need the Lindeberg condition 

\smallskip\noindent
(A0) $\ \, \limsup\nolimits_{n \to \infty}\ n^{-1}\, \sumnl_{k \le n}\,  E \,\{ d_k(\theta)^2\,  \one_{\,|d_k(\theta)| > \el\sqrt{n}}\} = 0 \quad (\theta \in \Theta,\ \el > 0).$

\begin{prop}\label{prop2}
Suppose there is some non-random function $\widetilde\gamma$ such that $\widetilde\gamma(\theta,\theta) >0,\, \theta \in \Theta$, and for every pair $\theta_1,\theta_2\in \Theta$ 
\beq\label{ascovt}
n^{-1}\ \sumnl_{k \le n} \widetilde d_k(\theta_1)\, \widetilde d_k(\theta_2) \equiv \widetilde\gamma_n(\theta_1,\theta_2)\ \longrightarrow_p\  \widetilde\gamma(\theta_1,\theta_2) \quad  \mbox{as} \ n \to \infty.
\eeq
Then under (A0), the finite-dimensional distributions of the process $\widetilde D_n$ converge to those of a mean zero Gaussian process $\widetilde Z$ with covariance $\widetilde \gamma$.
\end{prop}

The proposition suggests that for large $n$ the distribution of the test statistic $T(D_n)$ at the boundary, where $D_n = \widetilde D_n$, can be approximated by the distribution of the functional $T(\widetilde Z)$ on the paths of the Gaussian process $\widetilde Z$. Of course, convergence of the finite-dimensional distributions is insufficient for such a conclusion; tightness of the processes $\widetilde D_n$ in a suitable function space is required, too. Furthermore, the distribution of $T(\widetilde Z)$ generally is unknown and may be difficult to determine. And there still is the problem that the process $\widetilde D_n$ involves the (sum of the) conditional expectations $E\, [\, d_k(\theta)\, |\, \cF_{k-1}]$, which depend on the unknown probability $Q$ and would have to be estimated with sufficient accuracy. In view of these difficulties with the determination of proper critical values $\widetilde c_n$ we refer to the hypothetical test rejecting $H_-$ if $T(D_n)> \widetilde c_n$ as the ``fictitious test.''

For the sake of exposition we introduce further sub-hypotheses of $H_-$\footnote{In testing for forecast dominance the hypothesis $H_-$ is of our primary interest. The other hypotheses serve to develop the ideas step by step. In particular, the randomization test is not intended, nor apt, to test e.g.~$H_0^s$ vs.~$H_0$, or $H_-^s$ vs.~$H_-$.} besides $H_-^s$, namely the null-hypothesis $H_0$ of equal performance on average,
\beq\label{H0def}
\!\!\! H_0\, :\quad\, n^{-1/2}\, \sumnl_{k \le n} E\, [\, d_k(\theta)\, |\, \cF_{k-1}] \equiv M_n(\theta) = 0 \quad (\mbox{a.s.,} \ \theta\in \Theta),
\eeq
and the null-hypothesis $H_0^s$ of equal performance at every forecast instance,
\beq\label{H0sdef}
\!\!\! H_0^s\, :\quad\, E\, [\, d_k(\theta)\, |\, \cF_{k-1}] = 0 \quad (\mbox{a.s.,} \ \theta\in \Theta,\ k=1,\ldots,n).
\eeq
Evidently, $H_0,\, H_0^s$ may be regarded as the {\em boundaries} of the one-sided hypotheses $H_-,\, H_-^s$, respectively, representing different forms of equal predictive ability. Note that $H_0^s \subset H_-^s \subset H_- \subset H_-^w$, and $H_0 \subset H_-$; moreover, that a test that is valid for a given hypothesis remains valid when used as a test of a smaller sub-hypothesis; it may be invalid when used with a larger hypothesis comprising the first.

\section{Randomization tests}
\label{sec:rt}

\subsection{General idea}
\label{sec:gi}

Using a standard randomization device (c.f. Section \ref{sec:disc}), we let $\sigma_1,\, \sigma_2,\, \ldots$ be i.i.d.~such that $\sigma_k=\pm 1$ with probability 1/2 each, and define 
\beq\label{Drsdef}
D_n^*(\theta) = n^{-1/2}\, \sumnl_{k \le n}\, d_k(\theta) \sigma_k.
\eeq
We reject $H_-$ ``at level $\alpha$'' if $T(D_n) > c_n^*$ where $T$ is a monotone and convex test functional, and $c_n^*$ is determined such that $Pr^*[T(D_n^*) > c_n^*] \approx\alpha$. Here $Pr^*$ exclusively pertains to the random signs $\sigma_k$, the data $x_{k1}, x_{k2}, y_k$ being considered as fixed, non-random quantities. Henceforth we refer to this test as the {\em randomization test.} Its rationale is as follows. 
First, if $T$ is monotone, the argument at (\ref{alfacntrl}) {together with the inclusion $H_0 \subset H_-$ imply}
\beq\label{bdycontrol}
\supnl_{H_-} Pr\, [T(D_n) > c_n^*] = \supnl_{H_0} Pr\, [T(D_n) > c_n^*],
\eeq 
so that it suffices to control the error of the first kind across $H_0$ where there is no systematic difference between the two forecasts. Secondly, if there is no difference in predictive performance between forecasts 1 and 2, changing the labels (i.e., the sign of the $d_k$) should not affect the distribution of the test statistic. The quotes in ``at level $\alpha$'' shall underline that the test has (approximative) level $\alpha$ only {\em formally;} the actual test level might differ. 

While the randomization test has intuitive appeal and is easy to implement, its properties are less clear. For instance, calculating the critical value $c_n^*$ from the randomization distribution tacitly supposes that in the indifference case the distribution of the $\real^n$-valued process $\theta \mapsto (d_1(\theta),\ldots,d_n(\theta))$ is invariant under arbitrary sign changes in the $n$ components (same change for all $\theta$), which is an even stronger hypothesis than $H_0^s$.  This raises questions concerning the approximate range of validity of the test in asymptotic regimes, where fine distinctions between different hypotheses may become inessential. 

Initial answers will be obtained through (partially) heuristic reasoning, utilizing relationships between the randomization and the fictitious test. In Section \ref{sec:asy} these considerations are complemented by rigorous weak convergence results for the case of quantile and expectile forecasts, which validate the heuristics.

\subsection{Test validity: Heuristics, and a conjecture}
\label{sec:heu}

In part a) of the following proposition it is understood that chance enters in two ways: via the random signs $\sigma_k$, and via the statistical nature of the data triplets. In part b) we condition on the data, leaving the $\sigma_k$ as the sole source of randomness.

\begin{prop}\label{prop3}
a) Suppose there is a non-random function $\gamma$ such that $\gamma(\theta,\theta) >0,\, \theta \in \Theta$, and for every pair $\theta_1,\theta_2 \in \Theta$ 
\beq\label{ascov*}
n^{-1}\ \sumnl_{k \le n} d_k(\theta_1)\, d_k(\theta_2) \equiv \gamma_n(\theta_1,\theta_2)\ \longrightarrow_p \ \gamma(\theta_1,\theta_2).
\eeq
Then under assumption (A0) the finite-dimensional distributions of the process $D_n^*$ converge to those of a mean zero Gaussian process $Z$ with covariance $\gamma$.\\
b) The latter conclusion also holds under $Pr^*$ (i.e., conditionally on the data) provided that the stochastic convergence  (\ref{ascov*}) is replaced by the usual (deterministic) convergence, and the Lindeberg condition (A0) is satisfied without the expectation sign.
\end{prop}

\begin{rem}\label{rem1}
{\em Regarding part b), note that under $Pr^*$ the $d_k(\theta)$ are known non-random quantities, rendering the expectation sign void. On the other hand, if in (A0)\, $d_k(\theta)$ everywhere is replaced by $d_k^* (\theta)= d_k(\theta) \sigma_k$, and $E$ by the expectation $ E^*$ pertaining to the $\sigma_k$ only, then the resulting condition (A0$^*$) is a Lindeberg condition in the usual sense. Anyway, since $|d_k^* (\theta)| = |d_k(\theta)|$, there is no difference between the conditions with and without the expectation sign, and we need not distinguish (A0) and (A0$^*\!$).}
\end{rem}

We now address the question for which among the above hypotheses the randomization test is approximately valid, i.e., has size $\lessapprox \alpha$. The discussion builds on distributional approximations to be established later on and on an unproven conjecture. The argument still is instructive as it helps delineate the key problem. Let us re-emphasize that in regard to testing for forecast dominance our focus is on the one-sided hypothesis $H_-$; cf.~Footnote 5.

\smallskip\noindent
{\bf {\em Hypothesis $H_0^s$.}} \\
Under $H_0^s$ we have $d_k(\theta)=\widetilde d_k(\theta)$, hence $D_n =\widetilde D_n$ and $\gamma_n=\widetilde\gamma_n$. Consequently, the limit processes $Z$ and $\widetilde Z$ of $D_n^*$ and $\widetilde D_n$ are identical in distribution under $H_0^s$, and so are the limit distributions of any sufficiently regular test statistic (which need not be monotone, for instance). In particular, the critical values $c_n^*$ and $\widetilde c_n$ of the randomization and the fictitious test coincide asymptotically. Therefore, since the latter test is, for large $n$, approximatively valid for testing $H_0^s$ at level $\alpha$, then so is the former. The point is, of course, that the fictitious test is valid but infeasible, whereas the randomization test is straightforward to implement.

\smallskip\noindent
{\bf {\em Hypothesis $H_0$.}}\\
Under this hypothesis the above reasoning does not apply because the covariance functions $\gamma$ and $\widetilde \gamma$, hence the limit processes $Z$ and $\widetilde Z$, generally are different. Nevertheless, the randomization test remains approximatively valid for the hypothesis $H_0$ if the test functional $T$ fits with $H_0$. To substantiate this claim, let us begin by noting that $\gamma = \widetilde \gamma + \psi$ with a positive definite function $\psi$ given by the stochastic limit
\beq\label{psidef}
\psi(\theta_1,\theta_2) = p\,\mbox{-}\!\lim\,  n^{-1}\, \sumnl_{k \le n} E\, [\, d_k(\theta_1)\, |\, \cF_{k-1}] \ E\, [\, d_k(\theta_2)\, |\, \cF_{k-1}];
\eeq
cf.~Lemma \ref{lem3}. For the limit processes this means that $Z = \widetilde Z + W$ in distribution, where $W$ is an independent centered Gaussian process. \\
Now in the context of the hypotheses $H_0,\, H_0^s$ it is natural to consider test functionals $T$ that are {\em symmetric,} $T(-f) = T(f)$, and convex in $f$. In other words, the acceptance region $A = \{f:\, T(f) \le c\}$ is symmetric and convex. This allows us to control the error probability under $H_0$ by applying a celebrated inequality. A basic finite-dimensional version of the inequality is as follows.

\smallskip\noindent
{\bf Anderson's inequality} \cite[Corollary 3]{A}.
{\em Let $X, Y$ be independent centered $\real^d$-valued Gaussian random variables.  Let $g:\, \real^d \to \real$ be convex and symmetric, i.e.~$g(-x)=g(x)$ for every $x$. Then  $Pr[g(X) \le b] \ge Pr [g(X+Y) \le b]$ for every $b\in \real$.}
 
\smallskip\noindent
In our case $X$ and $Y$ correspond to $\widetilde Z$ and $W$ sampled discretely at $d$ points $\theta_j \in \Theta \subset \real$. Examples of functions $g$ corresponding to test functionals $T$ of interest are $g(x) = \max_i |x_i|$ and $g(x) = \sum_i |x_i|^p,\ p=1$ or $p=2$. 
As the sampling gets dense, one finds that under $H_0$ and for symmetric, convex test functionals $T$ one has
\beqa\label{knack}
Pr \left[T(D_n) > c_n^*\right] \!\!\! &  \stackrel{(1)}{=} & \!\!\! Pr \left[T(\widetilde D_n) > c_n^*\right] \approx Pr \left[T(\widetilde Z) > c_n^*\right] \\  && \!\!\!\!\!\!\!\!\!\!\! \stackrel{(2)}{\le}\,  Pr \left[T(Z) > c_n^*\right] \approx Pr^* \left[T(D_n^*) > c_n^*\right]  \approx \alpha .\nonumber 
\eeqa
Relation (1) holds because $D_n = \widetilde D_n$ under $H_0$, relation (2) by Anderson's inequality along with a standard approximation \cite[Proof of Corollary 4]{A}, and the three `$\approx$' signs hold by Propositions \ref{prop2} and \ref{prop3}, and by construction, respectively. Note that (\ref{knack}) in fact implies $\ \sup_{H_0} Pr \left[T(D_n) > c_n^*\right]\approx\alpha\, $ since $H_0 \supset H_0^s$ and the admissible error probability is fully exhausted on $H_0^s$.

\smallskip\noindent
{\bf {\em Hypotheses $H_-,\, H_-^s$.}}\\
The test functionals $T = T(f)$ appropriate for these hypotheses are convex and {\em monotone} in $f$. The latter property is incompatible with symmetry, which is an essential ingredient of Anderson's inequality. We nevertheless could argue similarly as above if there was a one-sided version of Anderson's inequality. The following would be most helpful.

\smallskip\noindent
{\bf A one-sided Anderson's inequality? -- Conjecture:}
{\em Let $X, Y$ be independent centered $\real^d$-valued Gaussian random variables.  Let $g:\, \real^d \to \real$ be convex and {\em monotone} in the sense that $g(x) \le g(y)$ whenever $x \le y$ (coordinatewise). Then there is a universal constant $\alpha_0 \in (0,1/2]$ (bold guess: $\alpha_0=1/2$) such that $Pr[g(X) \le b] \ge Pr [g(X+Y) \le b]$ whenever $Pr [g(X+Y) \le b] \ge 1-\alpha_0$.} 

\begin{rem}\label{rem2}
{\em An important consequence of this inequality is that {\em for test levels $\alpha \le \alpha_0$ the randomization test is (approximately) valid for testing the one-sided hypotheses} $H_-$, $H_-^s$. The argument is parallel to (\ref{knack}), with two modifications: first, since $M_n \le 0$ under $H_-, H_-^s$, instead of (1) we have an inequality, $Pr\,[T(D_n) > c_n^*] \le Pr\, [T(\widetilde D_n) > c_n^*]$,  by the monotonicity of $T$; secondly, relation (2) now follows from the one-sided Anderson inequality, again up to an approximation as in \cite[Proof of Corollary 4]{A}.}
\end{rem}

\begin{rem}\label{rem3}
{\em The one-sided Anderson's inequality is {\em not} needed for the approximative validity of the one-sided randomization test if we only consider (sequences of) probability measures $Q_n \in H_-$ that are {\em contiguous} \cite[p.~87]{vdV} to some sequence $P_n \in H_0^s$. This is because under $P_n$ we have $\gamma_n = \widetilde\gamma_n$, hence $\gamma_n \longrightarrow_p \gamma=\widetilde\gamma$, and by contiguity this convergence also takes place under $Q_n$; cf.~Propositions \ref{prop2}, \ref{prop3}. Therefore the distributions of the limit processes $Z$ and $\widetilde Z$ coincide, and the inequality (2) in (\ref{knack}) becomes an equality; while the `$=$' sign (1) there has to be replaced by `$\le$', again by the monotonicity of $T$. Thus in this case too, $Pr_{Q_n} \!\left[T(D_n) > c_n^*\right]\lessapprox \alpha$, as claimed. }
\end{rem}

\smallskip\noindent
{\bf Summary.}
{\em Asymptotically, the randomization test is an (approximatively) valid level-$\alpha$ test of the hypotheses $H_0,\, H_0^s$. If the conjecture is correct it is also valid for testing the hypotheses $H_-,\, H_-^s$.}   

\smallskip
It should be emphasized that the above discussion exclusively pertains to the control of the error of the first kind. Regarding power we only mention that the argument in \ref{rem3} may as well be applied to {\em alternatives} $Q_n \in H_+$ -- satisfying $M_{n,Q_n} \ge 0 \ Q_n$-a.s.~for all $\theta \in \Theta$; cf.~(\ref{H_def}) -- that are contiguous to a sequence $P_n \in H_0^s$. By monotonicity of $T$ one analogously gets $Pr_{Q_n} \!\left[T(D_n) > c_n^*\right]\gtrapprox \alpha$, i.e., {\em the randomization test is  unbiased against such alternatives.} For alternatives $Q_n$ not in $H_+$, where $\theta \mapsto M_{n,Q_n}(\theta)$ assumes positive as well as negative values, the power issue is subtle and will not be pursued systematically in this paper.

\subsection{Some comments on the conjectured inequality}
\label{sec:ineq}

In dimension $d=1$, the inequality is trivial. Convex, monotone acceptance regions then are intervals of the form $(-\infty,b]$, and if $X_i \sim \cN(0,v_i)\ (i=1,2)$ with $v_1  \le v_2$, then obviously $Pr[X_1 \le b] \ge Pr[X_2 \le b]$ if and only if $b \ge 0$, i.e., if and only if $Pr[X_2 \le b] \ge 1/2$. (This fits with the bold guess $\alpha_0=1/2$.)

For $d>1$, a small piece of evidence in favour of the conjecture can be given as follows. Let $A \subset \real^d$ be a convex acceptance region of the form $A\equiv A_{g,b} = \{x\in \real^d:\,  g(x) \le b\}$ for some convex function $g$ and $b \in \real$. 
(Monotonicity of $g$ is not required for the argument.) Denote by $G = \cN(0,V),\, G_+ = \cN(0,V_+)$ the distributions of the random variables $X$ and $X+Y$, respectively. The matrices $V,\, V_+$ and $V_+ - V$ are symmetric and (strictly) positive definite. Put $K = A \cap (-A)$, which intersection is convex and symmetric, and let $R$ denote the complement of the union $A \cup (-A)$. Then for any symmetric probability distribution $F$ on $\real^d$ we have
$
1 = F(A) + F(-A) - F(K) + F(R) = 2F(A) - F(K) + F(R)
$
or 
$
F(A) = (1+F(K) - F(R))/2.
$
For the moment being, suppose that $R$ is contained in the set $S$ where the density of $G_+$ exceeds the density of $G$. 
Then $G_+(R) \ge G(R)$, and an application of Anderson's inequality to the set $K$, $G_+(K) \le G(K)$, yields the desired conclusion,
$$
G_+(A) = (1+G_+(K) - G_+(R))/2 \le (1+G(K) - G(R))/2 = G(A).
$$
As for the possible inclusion $R \subset S$, note that in terms of the log densities 
$$
S = \{x:\, x'(V^{-1} - V_+^{-1})x   > \log \left(|V_+| /|V|\right)\}.
$$
Now $V_+  > V$ implies $V^{-1} > V_+^{-1}$ in the Loewner order \cite[Corollary 7.7.4(a)]{HJ}. Therefore $\Delta=V^{-1} - V_+^{-1}$ is positive definite, and noting that $L= \log \left(|V_+| /|V|\right)$ $> 0$ we find that $S$ is the complement of the ellipsoid $E = \{x:\, x'\Delta x \le L\}$. Since $A_{g,b} \uparrow \real^d$ as $b\uparrow \infty$ and $g$ is bounded on $E$, we have for all large enough $b$ that $S^c = E \subset A \subset R^c$, that is, $R\subset S$. But $b \to \infty$ iff the test level $\alpha \to 0$, so we have proved the following.

\begin{prop}\label{prop4}
If $A=A_{g,b}$ for some convex function $g$ and $b\in \real$, then there is $\alpha_0 \in (0,1)$ such that $Pr[g(X) \le b] \ge Pr [g(X+Y) \le b]$ whenever $Pr [g(X+Y) \le b] \ge 1-\alpha_0$.
\end{prop}

In our case, $X = (\widetilde Z(\theta_1),\ldots,\widetilde Z(\theta_d)),\ Y = (W(\theta_1),\ldots,W(\theta_d))$ with the $\theta_j$ becoming dense. Since the covariance function $\widetilde \gamma$ of $\widetilde Z$ generally is unknown, we have no control on the eigenvalues of $V^{-1} - V_+^{-1}$. Proposition \ref{prop4} thus does not guarantee that $\alpha_0$ stays bounded away from zero uniformly in the pair $V, V_+$ and all dimensions $d$, as it is necessary for the one-sided Anderson inequality. This uniformity is the core of the problem. 

A proof of the conjecture may require additional assumptions, e.g.~invariance of $g$ under coordinate permutations. (Generalizations involving other invariance conditions appear in \cite{DG,Mu}.) Relevant examples include the convex, monotone functions $g(x) = \max_i x_i,\, g(x) = \sum_i (x_i)_+^p\ (p \ge 1)$, which correspond to test functionals of major interest; cf.~(\ref{Tdef}). A proof for such a special case would already be most worthwhile. 

Thus far, our numerical experiments in the bivariate case $d=2$ and simulations with randomly generated covariance matrices for $d >2$ yielded no counterexample. Needless to say, this is irrelevant for the conjecture.

\section{Asymptotics for quantile and expectile forecasts}
\label{sec:asy} 

In principle, the developments so far apply to largely arbitrary functionals $\phi$ on the class $\cG$ of predictive distributions $G$ and related families of consistent scoring functions $S_\theta$. Hereafter we focus on functionals representing a quantile or an expectile. Given $\alpha \in (0,1)$, the $\alpha$-expectile of $G$ is defined as the unique solution $t$ to the equation $(1-\alpha)\int_{-\infty}^t (t-y)\, dG(y) = \alpha\int_t^\infty (y-t)\, dG(y)$ \cite{NP}; for $\alpha = 1/2$ one obtains the mean value functional. As usual, $q = \inf\{y:\, G(y) \ge \alpha\}$ is the (lower) $\alpha$-quantile of $G$, which here is identified with its right-continuous CDF. The median of $G$ obtains when $\alpha=1/2$.

As mentioned in Section \ref{sec:intro}, forecast dominance with respect to all consistent scoring functions is, for these functionals, equivalent to dominance with respect to a certain linearly indexed family of ``elementary'' scoring functions $S_\theta,$ such that any consistent scoring function can be written as 
\begin{equation}	
S(x,y) = \int S_\theta(x,y)\,  dM(\theta), \label{mixture}
\end{equation}
where $M$ is a nonnegative Borel measure on $\real$. Different functions $S(x,y)$ are characterized by different measures $M$. From \cite[Theorem 1]{EGJK}, the elementary scores for $\alpha$-quantiles are 
\beq\label{escq}
S_\theta(x,y) = \{\one_{y < x} - \alpha\} \, \{\one_{\theta < x} - \one_{\theta < y}\};
\eeq
setting $M$ in (\ref{mixture}) equal to the Lebesgue measure results in the popular `piecewise linear' score considered in quantile regression \cite{K1}. 
For $\alpha$-expectiles,
\beq\label{escx}
S_\theta(x,y) = |\one_{y < x} - \alpha\, | \, \{(y-\theta)_+  - (x-\theta)_+ - (y-x)\, \one_{\theta < x}\}; 
\eeq 
again, the most popular choice of $M$ is the Lebesgue measure, in which case $S(x,y)$ becomes the asymmetric squared error considered in \cite{NP}. 

The differences $d_k(\theta) = S_\theta(x_{k1},y_k) - S_\theta(x_{k2},y_k)$ of the elementary scores are distinguished by a particular property: $d_k(\theta)$ factorizes into the product of an {\em identification function} $I$ that depends only on the observation, times a difference of indicator functions depending only on the forecasts. Specifically, $d_k(\theta) = I(\theta,y_k)\{\one_{\theta < x_{k1}} - \one_{\theta < x_{k2}}\}$ \cite[Appendix A3]{EGJK}.
We shall utilize this fact to establish weak convergence results for quantile and expectile forecasts complementing those of Propositions \ref{prop2}, \ref{prop3} about finite-dimensional distributions.

\subsection{Conditional weak convergence of $D_n^*$}
\label{sec:condwc}

The purpose of this section is to establish the approximation $Pr^* [T(D_n^*) > c_n^*]$ $\approx Pr [T(Z) > c_n^*]$ figuring in the display (\ref{knack}) that is central to our argument. 
The asymptotics involves conditioning on the data $x_{k1},x_{k2},y_k$, so that the sign variables $\sigma_k$ form the only source of randomness. We thus avoid having to make assumptions about the stochastic structure of the data. 

Of basic importance are the second (cross-)moments of the process $D_n^*$,
\beqas
\gamma_n(\theta_1,\theta_2) \ceq E D_n^*(\theta_1) D_n^*(\theta_2) = n^{-1}\, \sumnl_{k \le n} d_k(\theta_1) d_k(\theta_2),\\ 
\rho_n(\theta_1,\theta_2)^2 \ceq E \left(D_n^*(\theta_2) - D_n^*(\theta_1)\right)^2 =  n^{-1}\, \sumnl_{k \le n}\, (d_k(\theta_2) - d_k(\theta_1))^2,
\eeqas
and the continuity moduli of the empirical distributions $G_n,\,F_{n1},\, F_{n2}$ of the observations $y_k$ and the forecasts $x_{k1},\, x_{k2}$, respectively.\\ Put $ m_k = |y_k -x_{k1}|\vee |y_k -x_{k2}|$, and let
\beqas
\mbox{for quantiles:} && 
H_n = G_n + F_{n1} + F_{n2}\, ,\\
\mbox{for expectiles:} && H_n = F_{n1} + F_{n2}\, .
\eeqas

\noindent
\textsc{Assumptions.}

\smallskip\noindent
(C1) \hspace*{1mm} (\ref{ascov*}) holds: there exists a function $\gamma$ such that $\gamma(\theta,\theta)>0,\,\theta \in \real$, and 
\beqas
n^{-1}\, \sumnl_{k \le n} d_k(\theta_1)\, d_k(\theta_2) \ceq \gamma_n(\theta_1,\theta_2)\, \longrightarrow \, \gamma(\theta_1,\theta_2) \quad (n \to \infty,\  \theta_1,\theta_2 \in \real). 
\eeqas

\noindent
(C2) \hspace*{1mm} There exist numbers $\kappa \in (0,1),\,B>0,\,  n_2 \ge 1$ and a sequence $\beta_n \to 0$ such that
$$
\supnl_{\,0 \le \theta_2 - \theta_1 \le r} \, H_n([\theta_1,\theta_2]) \le B\, (r\vee \beta_n)^\kappa, \quad r \in [0,1],\ n \ge n_2.
$$

\noindent
(C3)  \hspace*{17mm} $\sup_{n}\, n^{-1}\, \sumnl_{k \le n} m_k^4\, \equiv\, M \,<\, \infty$ \quad (only for expectiles). 

\medskip\noindent
(C4) \hspace*{1mm} There exist numbers $\nu> 0,\,A>0$, and  $n_1 \ge 1$ such that
$$
(F_{n1}+F_{n2})([-\theta,\theta]^c ) \le A\, \theta^{-\nu},\quad \theta \ge 1,\, n \ge n_1\, .
$$

\smallskip\noindent
{\bf {\em Discussion of the assumptions.}} 
In practice, the test functionals in (\ref{Tdef}) may be evaluated over a bounded interval $\Theta$ in the full $\theta$-domain $\real$ and the score difference processes restricted correspondingly. Of course, forecast dominance then holds only with respect to all mixtures $S = \int S_\theta\, dM(\theta)$ of scores $S_\theta$ with a measure $M$ supported by $\Theta$, which would often be considered as sufficient.\\
Assumption (C1) is a basic asymptotic stability requirement that would hold `in probability' under virtually any standard statistical model. However, as in part b) of Proposition \ref{prop3} 
there is no probability governing the data, hence no convergence in probability (see also Remark \ref{rem1}).\\
The uniform H\"older condition assumed in (C2) requires that the data $x_{k1},\, x_{k2},\, y_k$ are well dispersed and do not heavily cluster locally. The lower bound $r \ge \beta_n$ on the width of the increments is unavoidable because of the (asymptotically small) jumps of the empirical CDFs.\\
Assumption (C3) only matters for expectiles. To substantiate it, one may argue that reasonable forecasts should covary with the observations, which would limit the deflections of the quantities $m_k$. (C3) is stronger than boundedness on average of the $m_k^2$, which  appears as the minimal condition to impose. 
In return, it implies a Lindeberg type condition holding uniformly in $\theta$, 
\beq\label{LBC}
\lim\nolimits_{n \to \infty}\, \supnl_\theta\ n^{-1}\, \sumnl_{k \le n}\,  d_k(\theta)^2\,  \one_{\,|d_k(\theta)| > \el\sqrt{n}}\, =\, 0 \quad \mbox{for every $ \el > 0$.}
\eeq
Assumption (C4) restrains the large fluctuations of the forecasts $x_{k1},\, x_{k2}$ and allows us to control the tail behavior of the functions $\theta \mapsto E D_n^*(\theta)^2$.

\smallskip
Altogether, the assumptions appear weak as well as natural for the quantile and expectile functionals and for continuously distributed data. They only pertain to quantities computable from the data and do not presuppose any statistical model. On the other hand, if a probabilistic model is assumed, they hold with arbitrarily high probability in many of the customary settings. See the corresponding discussion in Section \ref{sec:addmat}, where (C2), (C4) are verified under conventional stationarity assumptions. 

Hereafter, $\ell_0^\infty$ denotes the space of all bounded measurable functions on $\real$ vanishing at infinity equipped with the sup-norm \cite{vdV}. The sample paths of $D_n^*$ are in $\ell_0^\infty$ since the elementary scoring functions of quantiles and expectiles are piecewise linear and vanish outside the smallest interval including all forecasts $x_{k\ell}$. In order to avoid problems related to the jumps of $D_n^*$, {\em we henceforth assume that the critical values $c_n^*$ are calculated using a continuous version $\bar D_n$ of $D_n^*$ obtained by linear interpolation of the $D_n^*$-values on the grid $\{j\beta_n:\, j \in \mathbb{Z}\}$, where $\beta_n$ is as in Assumption (C2). That is, the $c_n^*$ are supposed to satisfy} $Pr^*\,  [T(\bar D_n) > c_n^*] \approx \alpha$. Since the grid becomes arbitrarily fine as $n$ gets large this modification is immaterial; see Lemma \ref{lem2}.

\begin{thm}\label{thm1}
Under the assumptions (C1) to (C4) the processes $\{\bar D_n(\theta),\, \theta \in \real \}$ converge weakly in $\ell_0^\infty$ to a mean zero Gaussian process $\{Z(\theta),\, \theta \in \real\}$ with covariance function $\gamma$ and continuous sample paths.
\end{thm}

The theorem states weak convergence under $Pr^*$, i.e., conditionally on the data. A fortiori, this convergence holds unconditionally as long as the limit covariance function $\gamma$, hence the limit process $Z$, is unique. See Proposition \ref{prop3} for a related discussion.

As a consequence of the theorem, $T(\bar D_n)$ converges weakly in distribution to $T(Z)$ for any continuous functional on the space $\ell_0^\infty$. This covers the supremum statistic $T_\infty(f) = \sup_{\theta \in \real}\, f(\theta)$ as one special case of interest. Other examples such as the integral type functionals $T_p(f) = \int_\real f(\theta)_+^p\, d\theta$ ($p=1,2$) require a sharpening of assumption (C4) for the control of the tail masses.

\begin{cor}\label{cor1}
Assume (C1) to (C4). Then $T_1(\bar D_n)$ converges weakly in distribution to $T_1(Z)$ if the exponent $\nu$ in (C4) satisfies $\nu>2$ in the quantile, and $\nu>4$ in the expectile case. For the functional $T_2$ the corresponding conditions are $\nu>1$ in the quantile, and $\nu >2$ in the expectile case.
\end{cor}

The conditions on $\nu$ can be dropped if the integration in the functionals $T_1,\, T_2$ is restricted to a {\em finite} interval $\Theta \subset \real$ (cf.~(\ref{Tdef})). The simplification occurs because, in this case, the functionals are continuous on the space $\ell^\infty(\Theta)$ of all bounded measurable functions on $\Theta$ endowed with the sup-norm, and because weak convergence in $\ell_0^\infty$ implies weak convergence in $\ell^\infty(\Theta)$.

\subsection{Weak convergence of $\widetilde D_n$}
\label{sec:uncondwc}

Here the focus is on the approximation $Pr\, [T(\widetilde D_n) > c_n^*]  \approx Pr\, [T(\widetilde Z) > c_n^*]$ in (\ref{knack}), whose verification completes our proof of the validity of the randomization test (apart from the conjecture). The setting is unconditional, so '$Pr$' refers to some underlying probability measure governing the joint stochastic behavior of the data triplets. For simplicity we only deal with weak convergence on finite intervals $\Theta$, that is, in $\ell^\infty(\Theta)$.

The necessary distinction between the quantile and the expectile case is a bit tedious. We denote the sequentially conditioned versions of the empirical data distributions as 
$$
G_n^c(J) = \frac{1}{n}\, \sumnl_{k \le n} Pr[y_k \in J \mid \cF_{k-1}],\ F_{n\ell}^c(J) = \frac{1}{n}\, \sumnl_{k \le n} Pr[x_{k\ell} \in J \mid \cF_{k-1}], 
$$ 
$\ell = 1,2,\, J$ an interval, and put as earlier $H_n^c = G_n^c + F_{n1}^c + F_{n2}^c$ in the quantile, and $H_n^c = F_{n1}^c + F_{n2}^c$ in the expectile case. Note that $F_{n\ell}^c = F_{n\ell}$ in the common case of forecasts $x_{k\ell}$ that are $\cF_{k-1}$-measurable. 
The following assumptions are similar to those in the previous section, except that convergence is `in probabiliy' and expectations are being taken at the appropriate places. 
A justification of assumption (A2) is given in Section \ref{sec:addmat}.

\smallskip\noindent
\textsc{Assumptions.}

\smallskip\noindent
(A1) \hspace*{1mm} (\ref{ascovt}) holds: there exists a function $\widetilde \gamma$ such that $\widetilde\gamma(\theta,\theta)>0$ for all $\theta$, and 
\beqas
n^{-1}\, \sumnl_{k \le n} \widetilde d_k(\theta_1)\, \widetilde d_k(\theta_2) \ceq \widetilde \gamma_n(\theta_1,\theta_2)\, \longrightarrow_p \, \widetilde \gamma(\theta_1,\theta_2) \quad  (n \to \infty,\  \theta_1,\theta_2 \in \real).
\eeqas
\noindent
Given any $b>1$ there exists a number $p\ge 1$ such that (A2) and (A3) hold:

\smallskip\noindent
(A2) \hspace*{1mm} There are numbers $B>0,\, n_2 \ge 1$ and a sequence $\beta_n \to 0$ such \\
\hspace*{10mm} that for both $K_n = H_n$ and $K_n = H_n^c$
$$
\supnl_{\,0 \le \theta_2 - \theta_1 \le r} \, E K_n([\theta_1,\theta_2])^p \le B\, (r\vee \beta_n)^b, \quad r \in [0,1],\ n \ge n_2.
$$

\smallskip\noindent
(A3)  \hspace*{17mm} $\sup_{n}\, n^{-1}\, \sumnl_{k \le n} E m_k^{4p} \,<\, \infty$ \quad (only for expectiles).

\smallskip\noindent

As before we consider a continuous version of $\widetilde D_n$ obtained by linear interpolation on the grid $\{j\beta_n:\, j \in \mathbb{Z}\}$, $\beta_n$ as in Assumption (A2), which we denote as $\hat D_n$. Accordingly, the approximation to be established becomes $Pr\, [T(\hat D_n) > c_n^*]  \approx Pr\, [T(\widetilde Z) > c_n^*]$.

\begin{thm}\label{thm2}
Assume (A1) to (A3). For every bounded interval $\Theta \subset \real$ the processes $\{\hat D_n(\theta),\, \theta \in \real \}$ converge weakly in $\ell^\infty(\Theta)$ to a mean zero Gaussian process $\{\widetilde Z(\theta),\, \theta \in \real\}$ with covariance function $\widetilde \gamma$ and continuous sample paths.
\end{thm}

For bounded $\Theta$ the functionals $T_\infty,\, T_1,\, T_2$, and in fact all $T_p,\, 1 \le p \le \infty$ (see Equation \ref{Tdef}) are continuous on $\ell^\infty(\Theta)$. Therefore the desired approximation is immediate from the following.

\begin{cor}\label{cor2}
Assume (A1) to (A3). Then for any $1\le p \le \infty$ the random variable $T_p(\hat D_n)$ converges weakly in distribution to $T_p(Z)$.
\end{cor}

\section{Monte Carlo simulations}
\label{sec:mc}

Here we study the randomization test in finite sample scenarios involving mean (i.e., expectile) and quantile forecasts. The test statistics under examination are the positive part integrals $T_p(D_n) = \int D_n(\theta)_+^p\, d\theta,\ p =1,2$, considered as tests of the hypothesis $H_{-}^s$ saying that method 1 dominates method 2 at each time step. All simulation results are based on $1\,000$ Monte Carlo iterations.

\subsection{Mean forecasts}
\label{sec:mcm}
	
\begin{figure}[!h]
		\caption{Size and power of the randomization test for mean forecasts (top panel) and quantile forecasts (bottom panel).}
	\begin{tabular}{c}
		{\bf Mean forecasts} \\
		\includegraphics[width = \textwidth]{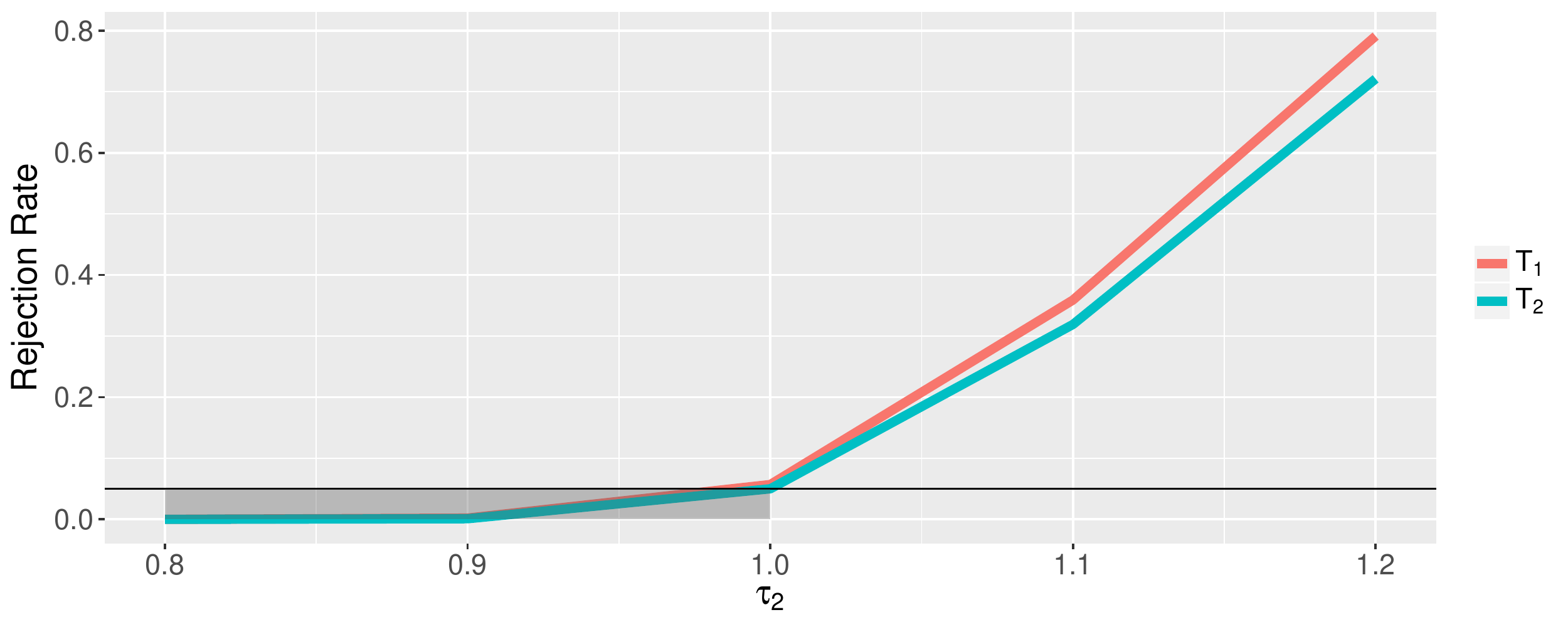} \\ 
		\\
		{\bf Quantile forecasts} \\
		\includegraphics[width = \textwidth]{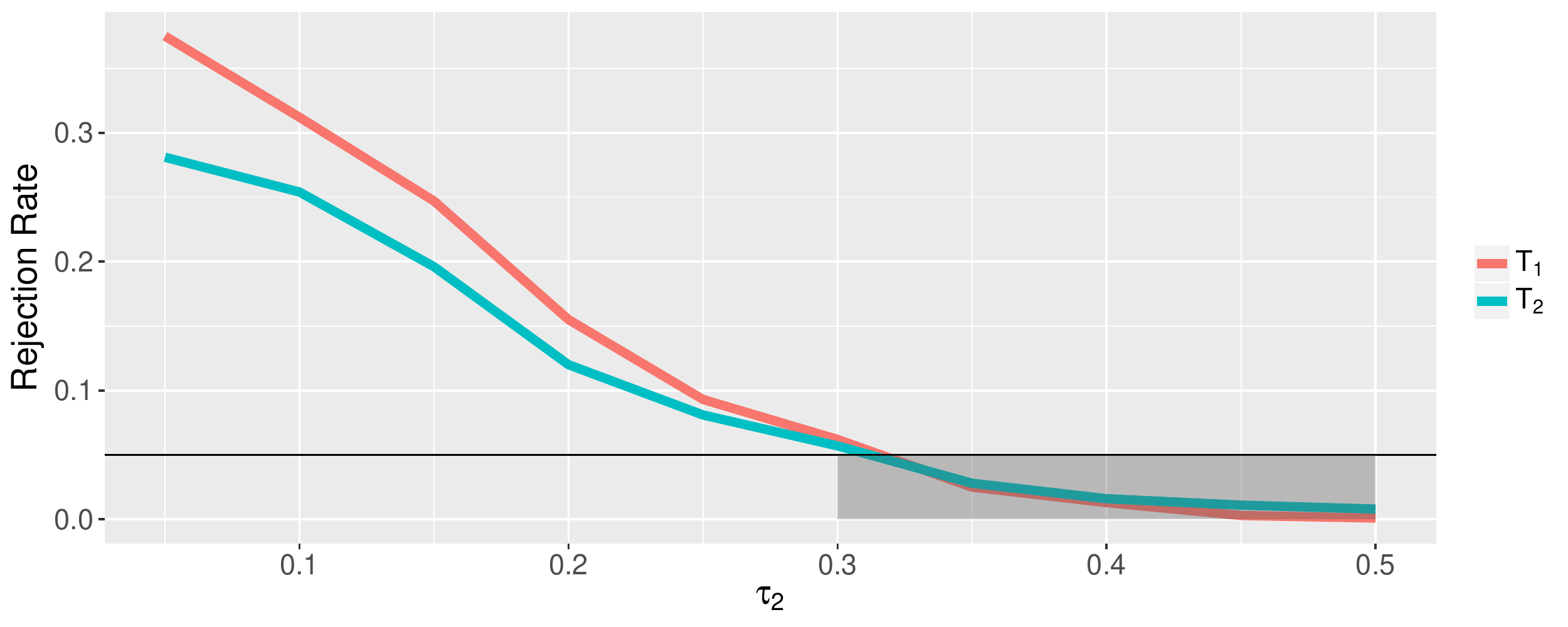} \\ [0.25cm]
	\end{tabular}
	{\em Rejections are at 5\% level, marked by horizontal line. In each panel, the dark gray area indicates the parameter range for which the hypothesis {$H_{-}^s$} is true, such that the rejection rate should be at most 5\%. Test statistics $T_1$ and $T_2$ are defined at (\ref{Tdef}). Results are based on $1\,000$ Monte Carlo iterations; within each iteration, the test is computed based on $1\,000$ simulated sign randomizations. See Sections \ref{sec:mcm} and \ref{sec:mcq} for further details. \label{fig:mc1}}
\end{figure}

We first present simulation results for the illustrative example from Section \ref{sec:nfd}. One of the variances is fixed,  $\tau_1 = 1$, while $\tau_2$ is varied. By Proposition \ref{prop1} the hypothesis $H_{-}^s$ is satisfied if $\tau_2 \le 1$, and violated otherwise.  We consider samples of $n = 200$ observations each, which in an economic context is empirically relevant for quarterly time series data focusing on the postwar period. The top panel of Figure \ref{fig:mc1} summarizes our results for the case where the regression parameter $a = 0.4$; similar results obtain for other values of $a$. The figure shows that the performance of the test is quite satisfactory: It comes close to its nominal level 5\% at the boundary of the hypothesis ($\tau_2 = 1$), and it is conservative in its interior ($\tau_2 < 1$), as predicted by the conjectured one-sided Anderson inequality. The part of the figure in which $\tau_2 > 1$ yields evidence on the power of the test. 
Naturally, we find that the power increases monotonically in $\tau_2$ (i.e., clearer violations of the hypothesis imply higher rejection rates). Furthermore, the functional $T_1$ has a slightly higher power than $T_2$.

\subsection{Quantile forecasts}
\label{sec:mcq}

We take the observations $y_k$ to follow an AR(1)-GARCH(1,1) process in the spirit of \cite{B}: 
\beqas
y_k \ceq 0.03 + 0.05\, y_{k-1} + s_k\varepsilon_k \\
s_k^2 \ceq 0.05 + 0.9\, s_{k-1}^2 + 0.05\, s_{k-1}^2 \varepsilon_{k-1}^2
\eeqas
with independent ``shocks'' $\varepsilon_k \sim \mathcal{N}(0,1)$. The parametrization follows \cite[Section 3]{Pt}, thereby intending to replicate the empirical features of daily stock returns. Forecasters are asked to state the $\alpha= 0.05$-quantile of the process, conditional on the information $\mathcal{F}_{k-1}$ available up to and including time $k-1$. 

To devise a simulation model for two imperfect forecasts, let us first conceive of an oracle. If the oracle knew the data generating mechanism and the initial values $s_0,\, \el_0,\, y_0$, she could successively compute $s_k$ from the observations $y_j,\, j <k$. Let $\cF_k$ denote the $\sigma$-algebra generated by the variables $s_0,\, \el_0,\, y_j,\, j \le k$. Then (assuming the regression parameters are known) $y_k \mid \cF_{k-1} \sim \cN(m_k,s_k^2)$, where we write $m_k = 0.03 + 0.05\, y_{k-1}$ for convenience. Thus for our oracle, the ideal quantile forecast would be the $\alpha$-quantile of the conditional distribution of $y_k$, namely $x_{k,ideal} = m_k + s_k z_\alpha$ where $z_\alpha = \Phi^{-1}(\alpha)$ is the ideal forecast in standard units.
This leads us to mimic lack of knowledge and forecast errors by assuming that the issued forecasts are of the form $x_{k\ell} = m_k + s_k z_{k\ell}\ (\ell =1,2)$ where the $z_{k\ell}$ are random perturbations of $z_\alpha$ that are independent among themselves and from all other variables. 
Specifically, we assume that the deflections from $z_\alpha$ are Gaussian in the log odds scale, 
\beq\label{qsim}
z_{k\ell} = \Phi^{-1}\!\left[\frac{1}{1+e^{-\beta-u_{k\ell}}}\right] \quad  (\beta = \log\frac{\alpha}{1-\alpha},\ u_{k\ell} \sim \cN(0,\tau_\ell^2),\ \ell = 1,2).
\eeq
Intuitively, forecast 1 should be better than forecast 2 if $\tau_1 < \tau_2$, since the deflections from the ideal forecast are then smaller for forecast 1. It can indeed be shown that $H_-^s$ holds if and only if $\tau_1 \le \tau_2$; cf.~end of Section \ref{sec:addmat}. 

In our simulations (bottom panel of Figure \ref{fig:mc1}), $\tau_1 = 0.3$ is fixed, and $\tau_2$ varies from $0.05$ to $0.5$. The quantile level is $\alpha = 0.05$, and the sample size is $2\,000$. Both choices are in line with the empirical case study in Section \ref{sec:csq}, where we analyze daily financial return data. Again, as in the previous example, the course of the power as a function of $\tau_2$ supports our claim that the randomization test is approximatively valid for testing $H_-^s$.

\section{Case studies}
\label{sec:cs}

\subsection{Mean forecasts}
\label{sec:csm}

For a practical application of the randomization test we consider the recession probability forecasts studied in \cite{RW}, using the updated data set analyzed by \cite[Section 4]{EGJK}. The data set covers $n = 186$ quarterly observations from 1968 to 2014, and two competing forecasting methods: Judgmental forecasts from a survey of professional forecasters (SPF), and forecasts from a simple statistical model (Probit). Both forecasts are one quarter ahead, and are out-of-sample.\footnote{The statistical model is re-estimated recursively at each forecast date in order to mimic a realistic forecast situation. The forecast data set is available within the \textsf{R} package \textsf{murphydiagram} \cite{JK}.} As shown in \cite[Figure 6]{EGJK}, the survey based forecasts attain better elementary expectile scores for most thresholds $\theta \in [0, 1]$. 
We specifically consider two test problems where either the survey forecast or the model based forecast dominates the respective other one under the maintained hypothesis. 

The top panel of Table 1 summarizes the results, which are based on $1\,000$ simulated sign randomizations. 
The hypothesis that Probit dominates SPF is rejected at the one percent significance level. By contrast, there is no evidence against the hypothesis that SPF dominates Probit. These results conform with those of \cite{EGJK} who consider informal (pointwise) confidence intervals. Remarkably, the randomization test here proves powerful enough to yield interpretable conclusions in a relatively small data set.

\begin{table}[!htbp]\label{tab:cs}
	\centering
		\caption{Randomization test results for recession probability forecasts (top panel) and quantile forecasts of stock returns (bottom panel).}
	\textbf{Mean forecasts (recession probabilities)} \\ [0.3cm]
	\begin{tabular}{p{4.5cm}cc}
		Hypothesis ($H_{-}$) & Test statistic & P-value \\ \toprule
		&& \\ [-0.2cm]
		SPF dominates Probit & $T_1$  & $0.987$ \\
		SPF dominates Probit & $T_2$ & $0.986$ \\
		Probit dominates SPF & $T_1$ & $0.010$ \\
		Probit dominates SPF & $T_2$ & $0.002$ \\ \bottomrule
	\end{tabular}\\ [0.4cm]
	\textbf{Quantile forecasts (stock returns, $\alpha = 0.05$)} \\ [0.3cm]
	\begin{tabular}{p{4.5cm}cc}
	 Hypothesis ($H_{-}$) & Test statistic & P-value \\ \toprule
		&& \\ [-0.2cm]
		$\text{QR}_{RV}$ dominates $\text{QR}_{ABS}$ & $T_1$  & $0.646$ \\
		$\text{QR}_{RV}$ dominates $\text{QR}_{ABS}$ & $T_2$ & $0.916$ \\
		$\text{QR}_{ABS}$ dominates $\text{QR}_{RV}$ & $T_1$ & $0.000$ \\
		$\text{QR}_{ABS}$ dominates $\text{QR}_{RV}$ & $T_2$ & $0.000$ \\ \bottomrule
	\end{tabular}\\ [0.3cm]
\end{table}

\subsection{Quantile forecasts}
\label{sec:csq}

In a second case study we consider quantile forecasts of daily returns $y_k$ of the Dow Jones Industrial Average (DJIA), using data that is freely available at \url{http://realized.oxford-man.ox.ac.uk/}. Quantiles at low levels $\alpha$ are commonly used as measures for financial risk, and are referred to as Value-at-Risk at level $\alpha$ (e.g. \cite{McNFE}, Sections 1 and 2). We specifically consider prediction of the five percent quantile of $y_k$, given information until the previous business day $k-1$. In a first specification, which we denote by $\text{QR}_{RV}$, the predicted quantile is given by
\begin{equation}
x_{k1} = \hat \beta_0 + \hat \beta_1 |RV_{k-1}|,\label{qrv}
\end{equation}
where $RV_{k-1}$ is the so-called realized volatility computed from intra-daily data (e.g. \cite{ABDL}). We obtain parameter estimates $\hat \beta_0, \hat \beta_1$ by quantile regression \cite{K1}, based on a rolling window of $2\,000$ observations.\footnote{The first rolling window ends on November 19, 2008; the last window ends on May 4, 2017. We use each rolling window to compute parameter estimates and form a forecast for the next business day. Our data set thus covers 1\,964 forecast/realization pairs that we use for evaluation. Our implementation of quantile regression is based on the function \textsf{rq} from the \textsf{R} package \textsf{quantreg} \cite{K2}.} Recent evidence \cite{ZB} suggests that the  specification in (\ref{qrv}) compares favorably to a number of more complicated alternatives. Our second specification ($\text{QR}_{ABS}$) is analogous to (\ref{qrv}), except that it employs the lagged absolute return $|y_{k-1}|$ in place of realized volatility $|RV_{k-1}|$. The two specifications are motivated by the fact that realized volatility and absolute returns proxy for the variability of financial returns, which is well known to fluctuate over time (cf.~Section \ref{sec:mcq}). Both measures should thus be informative about the quantiles of $y_k$, given $\mathcal{F}_{k-1}$.

The bottom panel of Table 1 presents the results of the comparison. We find no evidence against the hypothesis that $\text{QR}_{RV}$ dominates $\text{QR}_{ABS}$; however, we clearly reject the converse hypothesis that $\text{QR}_{ABS}$ dominates $\text{QR}_{RV}$. This suggests that intra-daily information encoded in realized volatility contains more predictive content than daily returns.\footnote{We obtain the same result when using the lagged squared return, $y_{k-1}^2$, rather than the lagged absolute return, $|y_{k-1}|$, in the second specification.} Similar conclusions were reached in \cite{ZKJF}. As in the first case study, the results are qualitatively robust across the two test statistics $T_1, T_2$. In summary, the Monte Carlo simulations and the case studies point to the potential usefulness of the proposed randomization test.

\section{Discussion}
\label{sec:disc}

We have studied randomization type tests of hypotheses implying that a quantile or expectile forecast is superior to a competitor, uniformly across all (or a large class of) consistent scoring functions. 
Variants of this topic recently have gained considerable interest, particularly in the econometrics literature. Tests of dominance relations in quantile and expectile forecasts are studied in \cite{Yen} using the bootstrap, while in \cite{ZKJF} the authors focus on the so-called expected shortfall functional, relying on a combination of pointwise tests and multiple testing corrections. These two papers are closest to the present work in that they base forecast comparisons on consistent scoring functions -- arguably the proper concept for this purpose \cite{G} -- and their mixture representation \cite{EGJK}. We next provide a more detailed discussion of related literature.

\textit{Dominance concept.} Our dominance concept aims at the direct comparison of forecasts on a purely empirical basis. This is distinct from more theoretically oriented concepts like nested or overlapping information sets \cite{HE,KrZi}, or the notion that forecast $x_1$ `encompasses' forecast $x_2$ \cite[Section 17.1]{ET}. In a nutshell, forecast encompassing asks whether a user who has access to both forecasts may safely ignore $x_2$, judging from some given criterion. By contrast, our notion of dominance asks whether a user who must choose between the forecasts would prefer $x_1$ over $x_2$, in regard to any admissible criterion.
Tests of {\em stochastic} dominance are considered in \cite{JCS,LMW,LSW}. Analogously to forecast dominance, stochastic dominance stands for superiority (of some procedure compared to another) that holds uniformly across a whole class of criteria. In stochastic dominance, this class comprises those functions of $e$, the residuals or forecast errors from some regression type model, that increase as $e$ moves away from zero, or else, are convex in $e$. In forecast dominance as here, the relevant criteria are the consistent scoring functions, or some subclass thereof. 

\textit{Modeling assumptions.} Our stance is to try to avoid assumptions about possible data generating mechanisms as far as possible, on the grounds given in the introduction. For a similar view see \cite{GW} and \cite[p.~1308 and Sect.~5]{JCS}, where the common stationarity assumption is weakened to distributional heterogeneity. Consequently our framework puts no restrictions on the type of forecasts or their connection with the observations, and it allows for great freedom regarding their dynamics. The lack of an explicit statistical model and the need to harmonize the underlying statistics with the surrogate external randomization lead us to state forecast dominance hypotheses in terms of {\em conditional} expectations expressing conditional predicitive ability, as proposed in \cite{GW}. The step-by-step character of the forecast scheme renders this an attractive, natural alternative to the familiar formulations using unconditional expectations; see Example \ref{xmpl} or the simulation model in Section \ref{sec:mcq}. We note, however, that our hypotheses differ from those of \cite{GW}: there the focus is on single-criterion two-sided (equality) testing, whereas in the present paper the focus is on simultaneous one-sided (inequality) testing. In the setup of \cite{GW}, simple least squares regressions yield attractive tests for the null-hypothesis of equal conditional predictive ability that are not available in our case.

\textit{Test procedures.} The idea of external randomization has long played an important role in statistics (see e.g. \cite[Chapter 4.4]{EH}, as well as \cite{CRom} on the related idea of permutation), and it
is central to our approach. Randomization for the task of forecast comparison was used in, e.g., \cite{DT,DM}, and it has been applied to ensure that a well-defined limiting distribution exists in the first place \cite{BC,CS,T}. Our test construction is similar in spirit to the `conditional p-value' approach of \cite{H} and the `wild bootstrap' proposed in \cite{M}. In the context of functional data, the use of maximum and integral type test statistics like ours is standard; for weighted versions cf.~\cite{LSW}. Since the limit distributions are unknown, the determination of critical values makes it necessary to resort to resampling methods, customarily various forms of (block \cite{Kue}) bootstrap as in \cite{JCS,LMW,LSW,Yen}. 
While this is often the method of choice, its application in the present one-sided, high-dimensional context is not without problems. These partly are due to the nonstandard asymptotics of the bootstrap-based tests resulting from degenerate limit processes; see e.g.~\cite{JCS,LSW,Yen}. Noteworthily at this point, the Gaussian limit processes in our setting are entirely regular. 

\textit{Test size control.} Intuition and classical test theory suggest that in order to control the error of the first kind, it might suffice to control it on the boundary of the hypothesis. Unfortunately, what constitutes the boundary is subtle, and a focus on least favorable cases is inadequate \cite{Hn}. For an extensive discussion of these points in a different framework (stochastic dominance) see Linton et al.~\cite{LSW}, who emphasize the importance, and difficulty, of a uniform control of the test size (see also \cite[p.~1320]{JCS}) and develop a sophisticated bootstrap procedure for this purpose (in the i.i.d.~case). Still, even there uniformity is achieved only if certain subsets of the hypothesis are excluded. We are actually not aware of any fully satisfactory result in this regard; neither is the issue clarified in the present paper. However, our approach suggests a potentially elegant solution at least: if the one-sided Anderson's inequality were true, our tests would be valid uniformly on $H_-$.\footnote{To be read as: uniformly on those parts of $H_-$ where the weak convergence in Theorems \ref{thm1}, \ref{thm2} holds uniformly. Clearly, such restrictions always apply when asymptotics is involved.} The discussion in Section \ref{sec:heu} elaborates the central role of our corresponding conjecture. A resolution of the issue, whether in the positive or in the negative, would certainly be of great interest. We may point out, however, that independently of the  final status of the conjecture, the randomization test does behave properly for probabilities that are contiguous to the strict null-hypothesis $H_0^s$; cf.~Remark \ref{rem3}.

\section{Proofs}
\label{sec:prfs}

{\bf {\em Proof of Proposition \ref{prop1}}.} We may refer to \cite[Example 2.3]{KrZi}, wherein our earlier proof is much simplified. To see the relevance of Theorem 2.1 of that paper for the elementary scores note that the mean value is an $\alpha=1/2$-expectile whence eq.~(\ref{escx}) assumes the form  $S_\theta = \varphi(y) - \varphi(x) - (y-x) \varphi'(x)$ with the convex function $\varphi(x) = \frac12\,(x-\theta)_+$. \done

\medskip\noindent
{\bf {\em Proof of Proposition \ref{prop2}.}} 
Given $\theta_1,\ldots,\theta_m \in \real$, $ c_1,\ldots,c_m \in \real$, put $X_n = \sum_j\, c_j \widetilde D_n(\theta_j)$ and $V = \sumnl_{i,j}\, c_i\, c_j\, \widetilde \gamma(\theta_i,\theta_j)$. It suffices to show that the distribution of $X_n$ converges to $\cN(0,V)$. We have $X_n = \sum_{k=1}^n X_{kn}$ where 
$$
X_{kn} = \sumnl_j\, c_j\,  n^{-1/2}\, \widetilde d_k(\theta_j), \quad k=1,\ldots,n.
$$
In order to apply \cite[Theorem 2.3]{McL} to the martingale difference array $\{X_{kn}\}$, we note at first that (\ref{ascovt}) implies
\beq\label{mcl1}
\sumnl_{k \le n}\, X_{kn}^2 = \sumnl_{i,j}\, c_i\, c_j\, n^{-1}\, \sumnl_{k \le n}\, \widetilde d_k(\theta_i)\,  \widetilde d_k(\theta_j) \longrightarrow_p \sumnl_{i,j}\, c_i\, c_j\, \widetilde \gamma(\theta_i,\theta_j) = V.
\eeq
Thus if 
\beq\label{mcl2}
E\left( \max\nolimits_{k\le n} X_{kn}^2\right) \to 0 \quad \mbox{as $n \to \infty$}
\eeq
holds, the two conditions (a), (b) of \cite[Theorem 2.3]{McL} are satisfied, and in view of (\ref{mcl1}) we are done. (We may asssume $V=1$.) Now
\beqas
\max\nolimits_{k\le n} X_{kn}^2 \cleq \frac{2m}{n}\, \sumnl_j c_j^2 \left\{\max\nolimits_{k\le n} d_k(\theta_j)^2 + \max\nolimits_{k\le n} E [\, d_k(\theta_j)^2 \, \mid\, \cF_{k-1}]\right\}
\eeqas
by Jensen's inequality, and since $m$ and the $c_j$ are fixed, it suffices to show that $n^{-1}$ times the expectation of the two maxima in curly brackets tends to zero for every $j$. Let $\el > 0$. For any $\theta$ we have 
\beqas
&& \!\!\! E \left(\max\nolimits_{k \le n}\,n^{-1} E \,[\, d_k(\theta)^2 \, \mid\,\cF_{k-1}] \right) \\ \cleq 
E \left( \max\nolimits_{k \le n}\, n^{-1} E \,[\, d_k(\theta)^2\, \one_{\,|d_k(\theta)| > \el\sqrt{n}} \, \mid\, \cF_{k-1}] \right)\nonumber\\
&& \! + E \left( \max\nolimits_{k \le n}\, n^{-1} E \,[\, d_k(\theta)^2\,  \one_{\,|d_k(\theta)| \le \el\sqrt{n}} \, \mid\, \cF_{k-1}] \right)\nonumber\\ \cleq
n^{-1}\, E \left( \sumnl_{k \le n}\,  E \,[\, d_k(\theta)^2\,  \one_{\,|d_k(\theta)| > \el\sqrt{n}} \, \mid\, \cF_{k-1}] \right) + \el^2 \nonumber\\ \ceq
n^{-1}\, \sumnl_{k \le n}\,  E \,\{ d_k(\theta)^2\,  \one_{\,|d_k(\theta)| > \el\sqrt{n}}\} + \el^2.
\eeqas
The same upper bound holds for $E \left(n^{-1} \max\nolimits_{k \le n} d_k(\theta)^2  \right)$. Since  $\el$ was arbitrary, (\ref{mcl2}) follows by assumption (C0). \done

\medskip\noindent
{\bf {\em Proof of Proposition \ref{prop3}.}} 
The proof follows the same lines as the proof of Proposition \ref{prop2}. It suffices to replace $\widetilde d_k(\theta)$ by $d_k(\theta) \sigma_k$, define $\cF_k$ as the $\sigma$-algebra generated by the random variables $\sigma_1,\ldots,\sigma_k$, and observe that  $|\sigma_k|=1$ and $E\, [\, d_k(\theta) \sigma_k\mid \cF_{k-1}]=0$. See Remark \ref{rem1}.
\done

\medskip
Toward the proofs of Theorem \ref{thm1} and Corollary \ref{cor1}, we assume throughout that (C1) to (C4) are fulfilled. We begin with some first consequences of the assumptions. Constants generally depend on whether they refer to quantiles or expectiles, which is indicated by subscripts $q,\, e$, respectively. Recall that expectations and probabilities are conditional, and refer to the random signs only. For simplicity we nevertheless use the ``un-starred'' symbols $E,\, Pr$.

\begin{lem}\label{lem1}
(i) 
There are constants $L_{q,e}$ such that 
\beq\label{l2bdd}
\supnl_n\, \int_{-\infty}^\infty E D_n^*(\theta)^2 d\theta \le L_{q,e}.
\eeq
(ii) 
There are constants $A_{q,e}$ and $\nu_q = \nu,\, \nu_e = \nu/2$ (cf.~(C4)) such that
\beq
ED_n^*(\theta)^2 \le A_{q,e}\, |\theta|^{-\nu_{q,e}},\ |\theta| \ge 1,\, n \ge n_1\, .\label{tvanish}
\eeq
(iii) There are constants $B_{q,e}$ and $\lambda_q = \kappa/2,\, \lambda_e = \kappa/4$ (cf.~(C2)) such that
\beq\label{holder}
\omega_n(r) = \supnl_{\,0 \le \theta_2 - \theta_1 \le r} \, \rho_n(\theta_1,\theta_2)  \le B_{q,e}\, (r\vee \beta_n)^{\lambda_{q,e}},\  r \in [0,1],\ n \ge n_2.
\eeq
\end{lem}

\noindent
{\em Proof.}
For generalized quantiles the individual score differences are of the form 
\beq\label{gscd}
d_k(\theta) = I(\theta,y_k)\,  \delta_k(\theta),\qquad \delta_k(\theta) = \one_{\theta < x_{k1}} - \one_{\theta < x_{k2}}
\eeq
where $I(\theta,y_k)$ is the respective identification function. Specifically for $\alpha$-quantiles, the identification function is $I(\theta,y) = \one_{y \le \theta} - \alpha$, whence $|d_k(\theta)| \le |\delta_k(\theta)|$. For $\alpha$-expectiles, $I(\theta,y) = (1-\alpha)\, (\theta- y)_+\, -  \alpha\, (y-\theta)_+$, whence 
\beq\label{dxbnd}
|d_k(\theta)| \le |y_k - \theta|\, |\delta_k(\theta)| \le \{|y_k -x_{k1}|\vee |y_k -x_{k2}|\}\, |\delta_k(\theta)| = m_k \, |\delta_k(\theta)|\, .
\eeq
The second inequality is easily seen to follow from the fact that $|\delta_k(\theta)|$ equals one if $\theta$ lies between $x_{k1}$ and $x_{k2}$, and is zero otherwise. 
This observation also shows that $\int \delta_k(\theta)^2 d\theta = |x_{k1}-x_{k2}| \le 2 m_k$, whence by H\"older's inequality and (C3) 
\beqas
\int E D_n^*(\theta)^2 d\theta \ceq n^{-1}\,\sumnl_{k\le n} \int d_k(\theta)^2 d\theta \le  n^{-1}\,\sumnl_{k\le n} 2(m_k)^{s+1}\le 2 M^{\frac{s+1}{4}}
\eeqas
where $s=0$ and $s=2$ for quantiles and expectiles, respectively, which is (i).
\\
Similarly, if $|\theta| \ge 1,\, n \ge n_1$, using (C4) we get for quantiles
$$
ED_n^*(\theta)^2 \le n^{-1}\, \sumnl_{k \le n}\,  |\delta_k(\theta)| \le (F_{n1}+F_{n2})([-|\theta|,|\theta|]^c ) \le A\, |\theta|^{-\nu}\, ,
$$
while for expectiles, (\ref{dxbnd}), (C3), and Cauchy-Schwarz give
\beqa
ED_n^*(\theta)^2 \ceq n^{-1}\, \sumnl_{k \le n}\,  d_k(\theta)^2 \le\left\{\left[n^{-1}\, \sumnl_{k \le n}\,  m_k^4\right] \left[n^{-1}\, \sumnl_{k \le n}\,  |\delta_k(\theta)|\right] \right\}^{1/2}\nonumber\\ \cleq \{MA\,  |\theta|^{-\nu}\}^{1/2}\, ,\label{meansq}
\eeqa
which settles (ii). As for the increments, let $\theta_1 < \theta_2$ and put $\delta_k(\theta_1,\theta_2) = \delta_k(\theta_2) - \delta_k(\theta_1)$. Writing 
$$
d_k(\theta_2) - d_k(\theta_1) = \{I(\theta_2,y_k) -I(\theta_1,y_k)\} \delta_k(\theta_i) + I(\theta_j,y_k)\delta_k(\theta_1,\theta_2)
$$
with either $i=1,j=2$ or $i=2,j=1$, whichever is more convenient, we get for $\alpha$-quantiles
\beq\label{iqest}
|d_k(\theta_2) - d_k(\theta_1)| \le 
\one_{\theta_1 <y_k \le \theta_2} + |\delta_k(\theta_1,\theta_2)|\, ,
\eeq
and for $\alpha$-expectiles
\beqa\label{ixest}
\!\!\!\!\!\!\!\!\!\! && \!\!\!|d_k(\theta_2) - d_k(\theta_1)| \\ \cleq  |\theta_2-\theta_1|(|\delta_k(\theta_1)| \vee |\delta_k(\theta_2)|) + \{|y_k -\theta_1|\wedge |y_k -\theta_2|\}\,  |\delta_k(\theta_1,\theta_2)| \nonumber\\ \cleq 
|\theta_2-\theta_1|(|\delta_k(\theta_1)| \vee |\delta_k(\theta_2)|) + \{|y_k -x_{k1}|\vee |y_k -x_{k2}|\}\,  |\delta_k(\theta_1,\theta_2)| .\nonumber
\eeqa
The last inequality may be verified similarly as at (\ref{dxbnd}) on observing that $|\delta_k(\theta_1,\theta_2)|=1$ if exactly one of $x_{k1}, x_{k2}$ lies in the interval $[\theta_1,\theta_2)$, and is zero otherwise. This observation also shows that 
\beq\label{delmean}
n^{-1}\,\sumnl_{k \le n}\,  |\delta_k(\theta_1,\theta_2)| \le F_{n1}([\theta_1,\theta_2)) + F_{n2}([\theta_1,\theta_2)).
\eeq
For quantiles we then get by (\ref{iqest})
\beqa
\rho_n(\theta_1,\theta_2)^2 \cleq 2\,  \{G_n((\theta_1,\theta_2]) + F_{n1}([\theta_1,\theta_2)) + F_{n2}([\theta_1,\theta_2)) \}\nonumber\\ \cleq  2 H_n([\theta_1,\theta_2]),\label{qrho2}
\eeqa
while for expectiles the estimates (\ref{ixest}), (\ref{delmean}) and Cauchy-Schwarz give similarly as at (\ref{meansq})
\beq\label{xrho2}
\rho_n(\theta_1,\theta_2)^2  \le 2\, (\theta_2-\theta_1)^2 + 2 \left\{M  H_n([\theta_1,\theta_2]\right\}^{1/2}.
\eeq
Assertion (iii) thus follows from (C2). 
\done

\begin{lem}\label{lem2}
Up to adjustments of the constants, the assertions of Lemma \ref{lem1} also hold for the interpolated processes $\bar D_n$, with the following improvement of (iii): 
\beq\label{iholder}
\supnl_{\,0 \le \theta_2 - \theta_1 \le r} \, E (\bar D_n(\theta_2) - \bar D_n(\theta_1))^2 \le C_{q,e}\, r^{\lambda_{q,e}} \quad (r \in [0,1],\ n\ge n_2)
\eeq
(i.e., with  $r^{\lambda_{q,e}}$ rather than $(r\vee \beta_n)^{\lambda_{q,e}}$). Furthermore,
\beqa\label{qmappr0}
\limnl_{n\to \infty}\ \supnl_\theta \, E (\bar D_n(\theta) - D_n^*(\theta))^2 = 0,\\
\limnl_{n\to \infty}\ E\! \int (\bar D_n(\theta) - D_n^*(\theta))^2 d\theta = 0.\label{qmappr1}
\eeqa
\end{lem}

\noindent
{\em Proof.} For convenience we intermediately write $\beta_n\equiv \el$. Given $\theta$, there is exactly one $\ell\in \mathbb{Z}$ and $w \in [0,1)$ such that $\theta = (1-w)\ell \el + w (\ell+1)\el$. By Jensen's inequality 
\beqas
&& \!\!\! E (\bar D_n(\theta) - D_n^*(\theta))^2 \\ \cleq  n^{-1}\,\sumnl_{k \le n}\left\{w\, [d_k((\ell+1)\el)-d_k(\theta)]^2 + (1-w)[d_k(\theta)- d_k(\ell\el)]^2\right\} \\ \cleq \omega_n(\el) \equiv \omega_n(\beta_n),
\eeqas
which proves (\ref{qmappr0}). Turning to (\ref{qmappr1}), let us write $\Delta_k(\theta)$ for the $k$-th term in the above sum. We first consider the quantile case. Recalling that $\ell$ and $w$ are uniquely determined by $\theta$ we get by (\ref{iqest}), 
$$
\Delta_k(\theta) \le 2\, \one_{\ell\el <y_k \le (\ell+1)\el} + 2\, |\delta_k(\ell\el,(\ell+1)\el)|\, .
$$

The right-hand side is always $\le 4$, and it vanishes except if both $\theta$ and any of $y_k,\, x_{k1}$, or $x_{k2}$ lie in the interval $[\ell\el,(\ell+1)\el)$. Thus for fixed $k$ there are at most $3$ intervals of length $\el$ on which the function $\theta\mapsto \Delta_k(\theta)$ is non-zero. Consequently, $\int \Delta_k(\theta)d\theta \le 12 \el$, so taking the average over $k$ settles the quantile case. A similar reasoning applies in the expectile case. By (\ref{ixest}),
$$
\Delta_k(\theta) \le 2\, \el^2\{|\delta_k(\ell\el)| \vee |\delta_k((\ell+1)\el)|\} + 2\, m_k^2\, |\delta_k(\ell\el,(\ell+1)\el)|\, .
$$
The term $|\delta_k(\ell\el)| \vee |\delta_k((\ell+1)\el)|\le 1$ is nonzero at most if $\theta \in [x_{k1}\wedge x_{k2}- \el,\,  x_{k1}\vee x_{k2}+ \el]$. Therefore
$$
\int \Delta_k(\theta)d\theta \le 2\, \el^2\, (|x_{k2}-x_{k1}|+2\el) + 4\,m_k^2\, \el \le 2\, \el^2\, (2m_k+2\el) + 4\, m_k^2\, \el,
$$
so averaging over $k$ and using (C3) gives (\ref{qmappr1}).\\
Straightforward estimates yield the uniform H\"older condition (\ref{iholder}) at first for points $\theta_1,\theta_2$ belonging to the same interval $[\ell\el,(\ell+1)\el)]$, then belonging to two adjacent intervals, finally for points with one or more such intervals in between, where we may apply (\ref{holder}). The analogs of assertions 1 and 2 in Lemma \ref{lem1} are obvious. 
\done

\medskip\noindent
{\bf {\em Proof of Theorem \ref{thm1}.}} 
Convergence of the finite-dimensional distributions being clear from Proposition \ref{prop3}, (\ref{LBC}), and (\ref{qmappr0}), we only need to prove (stochastic) asymptotic equicontinuity \cite{Po,vdV} and the uniform vanishing at infinity of the sample paths of $\bar D_n$. Without loss of generality we may assume $n \ge n_1\vee n_2$ (cf.~(C4), (C2)). Distinguishing between quantiles and expectiles is not necessary here, so we omit the subscripts $q,e$ in the quantities appearing in Lemma \ref{lem1} and \ref{lem2}. Moreover, by Lemma \ref{lem2} quantities initially referring to $D_n^*$ such as $\rho_n$ or $\omega_n$ may also be used with $\bar D_n$, with the same bounds.\\
Let $u>0$. For any set $T_0\subset\real$, let $N_n(u,T_0)$ denote the minimal cardinality of a subset $T \subset T_0$ such that $\min_{t \in T} \rho_n(\theta,t) \le u$ for every $\theta \in T_0$. Given $b>1$, pick $t_j \in [-b,b]$ equidistant with spacing $r = 2(u/B)^{1/\lambda}$. By (\ref{iholder}), the minimal $\rho_n$-distance of any $\theta \in [-b,b]$ to the resulting set $T$ is $\le \omega_n(r/2) \le u$, whence $N_n(u,[-b,b]) \le Kbu^{-1/\lambda}$. Here and subsequently we write $K$ for any independent finite constant, whose value may thus change from instance to instance.  \\
By (\ref{tvanish}) and Lemma \ref{lem2} there is $\nu >0$ such that
\beq\label{tailbnd}
\rho_n(\theta,b)^2 \le 2 \{E\bar D_n(\theta)^2 + E\bar D_n(b)^2\} \le (K b^{-\nu})^2,\quad \theta > b,
\eeq
and similarly for $\theta < -b$. Therefore, with $b = (K/u)^{1/\nu}$ we have $\min_{t \in T}\, \rho_n(\theta,t) \le u$ for every $\theta \in \real$ and thus 
\beq\label{coverbnd}
N_n(u,\real) \le K u^{-1/\nu-1/\lambda}, \quad u > 0.
\eeq
Let 
$$
\Omega_n(r) = \sup\,\{\, |\bar D_n(\theta_2)-\bar D_n(\theta_1)| :\,\rho_n(\theta_1,\theta_2) \le r,\,  \theta_1,\theta_2 \in \real\}, \ r >0.
$$ 
By \cite[Lemma 1.2]{Z} applied with $x_i(s) \equiv d_k(\theta)/\sqrt{n}$ and $p=1$ we have
$$
E\, \Omega_n(r) \le K \int_0^{r/4} (\log N_n(u,\real))^{1/2}\, du \le K \int_0^r (\log u^{-1/\tau-1/\lambda})^{1/2}\, du
$$
for all $r \in [0,1]$. Therefore $E\, \Omega_n(r_n) \to 0$ if $r_n \to 0$, which implies asymptotic equicontinuity on $\real$ with respect to the semi-metrics $\rho_n$.\\
There are two further consequences. First, we already know that for every $\eta >0$ there is  $b>0$ such that $E\bar D_n(b)^2 \le \eta^2$ and $\rho_n(b,\theta) \le \eta$ for every $\theta \in (b,\infty)$. Thus
$$
|\bar D_n(\theta)| \le |\bar D_n(b)| + |\bar D_n(\theta) - \bar D_n(b)| \le |\bar D_n(b)| + \Omega_n(\eta)
$$
and so  
$$
\supnl_{|\theta| >b,\, \theta \in \real}\, |\bar D_n(\theta)| = o_p(1) \quad \mbox{as}\quad  b \to \infty.
$$ 
Secondly, by (\ref{iholder})
\beqas
\widetilde\Omega_n(r) \ceq \sup\,\{\, |\bar D_n(\theta_2)-\bar D_n(\theta_1)| :\,|\theta_1 -\theta_2| \le r,\,  \theta_1,\theta_2 \in \real\} \\ \cleq \Omega_n(\omega_n(r)) \,\le\, \Omega_n( K r^\lambda),
\eeqas
whence $E\, \widetilde\Omega_n(r_n) \to 0$ if $r_n \to 0$, implying asymptotic equicontinuity also with respect to the standard metric. It follows that the processes $\bar D_n$ converge weakly in $\ell_0^\infty$ to the specified Gaussian process $Z$, which by the asymptotic equicontinuity can be assumed to have continuous sample paths.
\done

\medskip\noindent
{\bf {\em Proof of Corollary \ref{cor1}.}} By (\ref{tvanish}) and Lemma \ref{lem2}, $\int_{|\theta| > b}\, E \bar D_n(\theta)^2\, d\theta  \to 0$ as $b \to \infty$ under the given conditions. Consequently, $T_2(\bar D_n)$ equals $\int_{|\theta| \le b}\,  \bar D_n(\theta)_+^2\, d\theta$ up to the arbitrarily small contribution from the tails, and weak convergence follows by Theorem \ref{thm1}. The same argument, up to an application of Jensen's inequality, applies to $T_1$. \done

\medskip\noindent
{\bf {\em Proof of Theorem \ref{thm2}.}} 
To prove stochastic equicontinuity we use the classical Kolmogorov moment criterion. In view of the linear interpolation it suffices to show that there exist positive constants $\xi,K$, and $\eta>1$, such that for all $\theta_1,\theta_2$ in the grid $\{j\beta_n:\, j \in \mathbb{Z}\}$ one has
\beq\label{momcrit}
E\, |\widetilde D_n(\theta_2) - \widetilde D_n(\theta_1)|^\xi \le K |\theta_2-\theta_1|^\eta.
\eeq
Let such a pair $\theta_1,\theta_2$ be fixed. Since the partial sums 
$$
S_k = n^{-1/2}\,  \sumnl_{j \le k} (\widetilde d_j(\theta_2) - \widetilde d_j(\theta_1)), \quad k = 1,\ldots,n
$$
represent a martingale with respect to the filtration $\{\cF_k\}$, Burkholder's inequality \cite[Theorem 9]{Bu} gives 
\beqas
E\, |S_n|^{2p} \ceq E |\widetilde D_n(\theta_2) - \widetilde D_n(\theta_1)|^{2p} \\ \cleq
N_p\, E \left\{n^{-1} \sumnl_{k\le n} (\widetilde d_k(\theta_2) - \widetilde d_k(\theta_1))^2\right\}^p
\eeqas
for any $p \ge 1$, with a universal constant $N_p$. Now
\beqa
&& \!\!\! \widetilde \rho_n(\theta_1,\theta_2)^2 := n^{-1} \sumnl_k (\widetilde d_k(\theta_2) - \widetilde d_k(\theta_1))^2 \nonumber\\ \cleq
2 \rho_n(\theta_1,\theta_2)^2 +2 n^{-1} \sumnl_k \{E\,[\,|d_k(\theta_2) - d_k(\theta_1)|\mid \cF_{k-1}]\}^2\label{wtrhoest}
\eeqa
which may be further estimated as in the proof of Lemma \ref{lem1}.  
We first consider the expectile case. Putting $\delta_k = \delta_k(\theta_1,\theta_2)$ we get from (\ref{ixest}) that the last term is bounded by a constant times the sum of $(\theta_2-\theta_1)^2$ plus the term 
\beqas
&& \!\!\! n^{-1} \sumnl_k \{E [\,m_k\, |\delta_k|\mid \cF_{k-1}]\}^2 \le
n^{-1} \sumnl_k E [\,m_k^2 \mid \cF_{k-1}]\, E [\,\delta_k^2\mid \cF_{k-1}] \\ \cleq
\left\{n^{-1} \sumnl_k (E [\,m_k^2 \mid \cF_{k-1}])^2\right\}^{1/2}\left\{n^{-1} \sumnl_k E [\,|\delta_k|\mid \cF_{k-1}]\right\}^{1/2}
\\ \cleq
\left\{n^{-1} \sumnl_k E [\,m_k^4\mid \cF_{k-1}]\right\}^{1/2}\left\{H_n^c([\theta_1,\theta_2])\right\}^{1/2}.
\eeqas
It follows that
\beqas
E\, \widetilde \rho_n(\theta_1,\theta_2)^{2p} \cleq 
K \left[ |\theta_2-\theta_1|^{2p} + E \left\{\big(n^{-1} \sumnl_k m_k^4\big)\, H_n([\theta_1,\theta_2]\right\}^{p/2}\right.\\
&& \left.\qquad\ +\, E \left\{\big(n^{-1} \sumnl_k E [\,m_k^4\mid \cF_{k-1}]\big)\, H_n^c([\theta_1,\theta_2]\right\}^{p/2} \right] \\ \cleq
K \left[ |\theta_2-\theta_1|^{2p} + \left\{E\, \big(n^{-1} \sumnl_k m_k^4\big)^p\ E\,H_n([\theta_1,\theta_2]^p\, \right\}^{1/2}\right.\\
&& \left.\qquad\ +\, \left\{E\, \big(n^{-1} \sumnl_k E [\,m_k^4\mid \cF_{k-1}]\big)^p\ E\, H_n^c([\theta_1,\theta_2]^p\right\}^{1/2} \right]
\\ \cleq
K \left[ |\theta_2-\theta_1|^{2p} + \sqrt{M_{4p}}\, \big(E\,H_n([\theta_1,\theta_2]^p+ E\, H_n^c([\theta_1,\theta_2]^p\big)^{1/2} \right]
\eeqas
where $M_{4p} = n^{-1} \sumnl_k E\, m_k^{4p}$. So given $\eta>1$, putting $b = 2\eta$ we may choose $p\ge 1$ in (A2), (A3) such that (\ref{momcrit}) is satisfied with $\xi = 2p$. This settles the expectile case. The quantile case can be dealt with similarly starting from (\ref{wtrhoest}). Given $\eta>1$ one puts $b =\eta$ and uses (\ref{iqest}), (\ref{qrho2}), then (A2), (A3).
\done

\section{Additional material}
\label{sec:addmat}

\begin{lem}\label{lem3}
Let $h_k(\theta) = E\, [\, d_k(\theta)\, |\, \cF_{k-1}]$. Under the conditions (\ref{ascovt}), (\ref{ascov*}), and (C3) we have
\beq
\gamma = \widetilde \gamma + \psi\quad where \quad  
\psi(\theta_1,\theta_2) = p\,\mbox{-}\!\lim\,  n^{-1}\, \sumnl_{k \le n} h_k(\theta_1) h_k(\theta_2).
\eeq
\end{lem}

\noindent
{\em Proof.} By (\ref{ascovt}), (\ref{ascov*}), and
\beqas
\gamma_n(\theta_1,\theta_2) \ceq n^{-1}\, \sumnl_{k \le n}\{\widetilde d_k(\theta_1) + h_k(\theta_1)\}\,\{  \widetilde d_k(\theta_2) + h_k(\theta_2)\}\\ \ceq
\widetilde\gamma_n(\theta_1,\theta_2) + n^{-1}\, \sumnl_{k \le n} h_k(\theta_1) h_k(\theta_2)\\
&&\!\!\! + \,  n^{-1}\, \sumnl_{k \le n} h_k(\theta_2)\widetilde d_k(\theta_1) + n^{-1}\, \sumnl_{k \le n} h_k(\theta_1)\widetilde d_k(\theta_2)
\eeqas
it suffices to show that e.g.~the last term, to be denoted $R_n$, tends to zero in quadratic mean. But $E R_n = 0$ because $E\, [\, \widetilde d_k(\theta_2)\, |\, \cF_{k-1}] = 0$ and $h_k(\theta_1)$ is $\cF_{k-1}$-measurable. Similarly, $E R_n^2 \to 0$: the off-diagonal terms in the double sum vanish, and by Jensen's and Cauchy's inequalities and (C3) the sum of the diagonal terms is $O(n)$. \done

\medskip\noindent
{\bf {\em Justification of (C2), (C4).} } 
We will show that the conditions (C2), (C4) are satisfied with probability arbitrarily close to one under common probability models for the data. Possible dependencies within the triplets $(x_{k1},x_{k2},y_k)$ do not matter because (C2), (C4) effectively pertain to the marginal CDFs $F_{n1},F_{n2},G_n$ only. However, it is natural in the prediction setting to allow for serial dependence. 
Specifically, suppose that the predictions $x_{k1},x_{k2}$ and the observations $y_k$ each form a strictly stationary sequence, defined for all $k \in \mathbb{Z}$.

\smallskip\noindent
CONDITION (C2). To verify (C2) for the empirical CDF $G_n$ of the observations, e.g., we may invoke an estimate  by W.~B.~Wu applying to certain causal processes of the form $y_k = J(\cdots,\el_{k-1},\el_k)$, where $J$ is measurable and the $\el_k, \, k \in \mathbb{Z}$ are i.i.d.~random variables. As an immediate consequence of \cite[Theorem 2]{Wu} one has, under the conditions given there, that 
$$
E \left[\,\supnl_{\,0 \le t-s \le r}\, (G_n(t) - G_n(s))^2\right] = O(n^{-1}\, r^{1-2/q}),\qquad r \ge \beta_n,
$$
where $2 < q < 4$ and $\beta_n$ is a sequence tending to zero sufficiently slowly (not faster than $\,n^{-1}\,(\log n)^{2q/(q-2)}\, $).
Markov's inequality then gives
$$
Pr\left[\,\supnl_{\,0 \le t-s \le r}\, |G_n(t) - G_n(s)| \ge  K r^\kappa\right] = O(n^{-1}\, r^{1-2/q-2\kappa}),  \quad r \ge \beta_n,\, K > 0.
$$
Putting $r = 2^{-\ell}$ and summing over $\ell= 0,1,\ldots$ we find that for any positive $\kappa < 1/2 - 1/q$ we have with probability $1 - O(n^{-1})$ that 
$$
\supnl_{\,0 \le t-s \le 2^{-\ell}}\, |G_n(t) - G_n(s)| \le K (2^{-\ell})^\kappa \quad\mbox{for every $\ell$ with } 2^{-\ell} \ge \beta_n.
$$
Thus (C2) holds with probability tending to one for the processes in question.

For an alternative justification let us consider the more common case where the $y_k$ form a strong ($\alpha$-)mixing sequence. In a first step we apply covariance inequalities due to E.~Rio \cite{Ri1} yielding an estimate of the variance of the increments of $G_n$. Specifically, suppose that the mixing coefficients $\alpha_n$ decay as $n^{-\varrho}$ for some $\varrho > 1$. Let $G$ denote the common CDF of the $y_k$. Then as a consequence of \cite[Theorem 1.2]{Ri1} we get that for some finite constant $K$
\beq\label{Rio1}
\mathrm{Var}(G_n(t) - G_n(s)) \le  K n^{-1} (G(t) - G(s))^{1-1/\varrho}, \quad s < t,\, n \ge 1.
\eeq
{\em Proof.} (Cf.~\cite[pp.~590,591]{Ri1}.) We have $G_n(t) - G_n(s) = n^{-1}\sumnl_{1 \le k\le n}\,  \xi_k,\ \xi_k = \one_{s < y_k \le t}$. The common `quantile function' $Q_\xi$ of the $\xi_k$ is readily seen to be given by $Q_\xi(u) = 1$ if $u < Pr[s < y_k \le t] =  G(t) - G(s) \equiv \Delta$, and $=0$ otherwise. Since $\alpha^{-1}(u) = \sumnl_k\, \one_{\alpha_k >u} = O(u^{-1/\varrho})$ it follows that the term
$$
\int_0^1 \alpha^{-1}(u) Q_\xi(2u)^2\, du \le K \int_0^\Delta u^{-1/\varrho}\, du =  O(\Delta^{1-1/\varrho}), 
$$
whence (\ref{Rio1}) follows. \done

\noindent
Further by \cite[Theorem 1.2]{Ri1}, the limits of the sequences $n\,\mathrm{Var}(G_n(t) - G_n(s))$ and $n\,\mathrm{Var}(G_n(t))$, hence also of $n\,\mathrm{Cov}(G_n(s),G_n(t))$, exist for all $s,t$. The latter, for instance, is given by the absolutely convergent sum
\beq\label{Lamcov}
\limnl_{n\to \infty}\, n\,\mathrm{Cov}(G_n(s),G_n(t)) = \sumnl_{k \in \mathbb{Z}} \mathrm{Cov}(\one_{y_0 \le s}, \one_{y_k \le t}) \equiv \Lambda(s,t).
\eeq
Analogous expressions hold for the other limits.\\
We next appeal to the weak convergence of the processes $n^{1/2} (G_n(t) - G(t)),\, t \in \real$ to the mean zero Gaussian process $V(t),\, t \in \real$ with covariance function $\Lambda$ from (\ref{Lamcov}), which follows from a stronger (almost sure) approximation result cited below. 
By the convergence of moments, the increments of $V$ satisfy the analog of (\ref{Rio1}),
$$
E (V(t) - V(s))^2 \le K (G(t) - G(s))^{1-1/\varrho}.
$$
Since $V$ is Gaussian, an application of the well-known Garsia-Rodemich-Rumsey Lemma (see e.g.~\cite{ASVY}) along with an intermediate time change implies that there exists a positive constant $K$ such that for any $\delta < (1-1/\varrho)/2\, $ the process $V$ satisfies with probability one a pathwise H\"older condition of the form 
\beq\label{GHold}
|V(t) - V(s)| \le K (G(t)- G(s))^\delta, \quad s < t \quad  \mbox{(almost surely, `a.s.')}.
\eeq
Now suppose that the CDF $G$ is uniformly H\"older continuous with index $\kappa \in (0,1]$. Then by (\ref{GHold}), $V$ also fulfils, for any $\eta < \kappa(1-1/\varrho)/2\, $, 
$$
|V(t) - V(s)| \le K (t-s)^\eta, \quad s < t \quad  \mbox{(a.s.)}.
$$

\noindent
In order to transfer this pathwise H\"older condition to the processes $G_n$ we may apply a ``Hungarian type'' strong approximation result for the empirical process of a stationary sequence. As a consequence of \cite[Theorem]{Y} or \cite[Theorem 2.1]{DMR}, there exists a sequence of Gaussian processes $V_n(t),\, t \in \real$, all copies of $V$, all defined on a common probability space carrying also the $y_k$, such that 
\beq\label{strapp}
\supnl_{t\in \real}\, |G_n(t)-G(t) - n^{-1/2} V_n(t)| = o(n^{-1/2}) \qquad \mbox{(a.s.)}.
\eeq
It follows that a.s.~for all $s,t$ with $|s-t| \le r$
\beqas
|G_n(t) - G_n(s)| \cleq |G(t) - G(s)| + n^{-1/2}|V_n(t) - V_n(s)| + o(n^{-1/2})\\ \ceq
O(r^\kappa) +  O(n^{-1/2}\, r^{\eta}) + o(n^{-1/2}).
\eeqas
This reduces to $O(r^\kappa)$ if we set $\beta_n = n^{-1/(2\kappa)}$, since then $n^{-1/2} = \beta_n^\kappa \le r^{\kappa}$ for $r \ge \beta_n$. The empirical processes $F_{n1},\, F_{n2}$ can be treated analogously. Thus we have shown that assumption (C2) is fulfilled with probability one under the indicated conditions, namely (sufficiently) strong mixing of the $y_k$ and $x_{k1},\, x_{k2}$, and H\"older continuity of their respective marginal CDFs.  \done

\smallskip\noindent
CONDITION (C4). 
Suppose that the predictions $x_{k1}$ form a strongly mixing sequence with the common marginal CDF $F_1$.  
Suppose, furthermore, that we have a strong approximation of the empirical processes $F_{n1}$ as in the previously discussed case. Using the same notation $V_n$ for the approximating Gaussian processes, we then have as in (\ref{strapp})
$$
\supnl_{\theta \in \real}\, |F_{n1}(\theta) - F_1(\theta) - n^{-1/2}\, V_n(\theta)| = o(n^{-1/2}) \quad \mbox{(a.s.)}. 
$$
Arguing as from (\ref{Rio1}) to (\ref{GHold}) we get $|V_n(-\theta)| \le K F_1(-\theta)^\delta$ for all $\theta\ge 0$ (a.s.), where again $\delta$ may be any positive number less than $(1 - 1/\varrho)/2$ and $\varrho > 1$ has the same meaning for the $x_{k1}$ as it had for the $y_k$. \\
We now assume that $\int |x|^q dF_1(x) < \infty$ for some $q >0$. Then 
\beqas
F_{n1}(-\theta) \cleq F_1(-\theta) + n^{-1/2}\, |V_n(-\theta)| + o(n^{-1/2}) \\ \ceq O(\theta^{-q}) + O(n^{-1/2}\, \theta^{-q\delta}) + o(n^{-1/2}), \quad \theta \ge 1 \quad \mbox{(a.s.)}
\eeqas
For $1 \le \theta \le n^{1/q}$, hence $n^{-1/2} \le \theta^{-q/2}$, we have 
$$
F_{n1}(-\theta) = O(\theta^{-q}) + O(\theta^{-q/2 - q\delta})+ o(\theta^{-q/2}),
$$
which certainly is $o(\theta^{-\nu})$ if we set $\nu = q/2$. An analogous estimate for the right tail gives 
$F_{n1}([-\theta,\theta]^c) = o(\theta^{-\nu})$ uniformly in the range $1\le \theta \le n^{1/q}$ (a.s.). On the other hand, if $\theta > n^{1/q}$ then
$$
Pr\, [F_{n1}([-\theta,\theta]^c) > 0] = Pr\, [\maxnl_{k \le n}\, |x_{k1}| > \theta] \le n\theta^{-q}\int_{|x| > n^\frac{1}{q}} |x|^q dF_1(x),
$$
which is $o(1)$ as $n\to\infty$. It follows that with probability tending to one we have $F_{n1}([-\theta,\theta]^c) = o(\theta^{-\nu})$ for all $\theta \ge 1$. Thus under the indicated assumptions, the tail condition (C4) is fulfilled with arbitrarily high probability if the marginal CDFs $F_1,F_2$ have a finite absolute moment of the order $q = 2\nu$.   \done

\medskip\noindent
{\bf {\em Justification of assumption (A2).} } 
In view of general Poisson approximation results for frequencies of rare events (e.g., \cite{CR}) we may expect that under broad conditions allowing for dependent observations the number $N$ of data falling into a small interval $J$ of length $r$ is roughly Poisson distributed with parameter of the order $n r^\kappa$, where $\kappa \in (0,1]$ characterizes the (maximal) clustering of the data points. 
The $p$-th moment $(p\ge 1)$ of the Poisson distribution with mean $m$ is $O(m+m^p)$ uniformly in $m$. Thus if $nr^\kappa \ge 1$ we may expect that $E H_n(J)^p = O(n^{-p} (nr^\kappa)^p) = O(r^{p\kappa})$, while for $nr^\kappa \le 1$ we should have $E H_n(J)^p = O(n^{-p} nr^\kappa) = O(n^{1-p} r^\kappa)$. 
Now given $b>1$, choose $p > b/\kappa$ and put $\beta_n = n^{-(p-1)/(b-\kappa)}$. Noting that $\beta_n \le n^{-1/\kappa}$, we find that in case $nr^\kappa \ge 1$ we have both $E H_n(J)^p = O(r^{p\kappa}) = O(r^b)$ and $r \ge \beta_n$. In case $nr^\kappa \le 1$ we have $E H_n(J)^p = O(n^{1-p} r^\kappa)$, which is $O(\beta_n^{b-\kappa} r^\kappa) = O(r^{b})$ if $r \ge \beta_n$. 
It follows that $E H_n(J)^p = O(r^b)$ whenever $r \ge \beta_n$, so that (A2) indeed would hold under quite general conditions. \done

\medskip\noindent
{\bf {\em Analysis of the quantile forecast example (Section \ref{sec:mcq}).} }
The difference of the elementary quantile scores is
$$
d_k(\theta) = S_\theta(x_{k1},y_k) - S_\theta(x_{k2},y_k) = (\one_{y_k \le \theta} - \alpha)\, \{\one_{\theta < x_{k1}} - \one_{\theta < x_{k2}}\}.
$$
Taking our assumptions into account and passing to standard units on writing $t_k = (\theta - m_k)/s_k$ (and $z_{k\ell} = (x_{k\ell} - m_k)/s_k$), we get
\beq\label{expsd}
h_k(\theta) := E\, [\, d_k(\theta)\mid \cF_{k-1}] = (\Phi(t_k) - \alpha)\, \{Pr[t_k < z_{k1}] - Pr[t_k < z_{k2}]\}
\eeq
where $Pr$ refers to the $z_{k\ell}$ (resp.~$u_{k\ell}$), everything else being considered as nonrandom (given $\cF_{k-1}$). We henceforth omit the index $k$ and use the abbreviation $h_k(\theta) \equiv h$. \\
Since $t < z_\ell$ iff $\log[\Phi(t)/(1-\Phi(t))] - \log[\alpha/(1-\alpha)] < u_\ell$, we have with $\lambda(p) = \log[p/(1-p)]$ that
$$
Pr[t < z_\ell] = 1-\Phi\left([\lambda(\Phi(t))-\lambda(\alpha)]/\tau_\ell\right) = \Phi\left([\lambda(\alpha)- \lambda(\Phi(t))]/\tau_\ell\right).
$$
Suppose at first that $\Phi(t) - \alpha <0$. Then $\lambda(\alpha)- \lambda(\Phi(t)) >0$, so $\tau_1 < \tau_2$ implies 
$$
Pr[t < z_1] = \Phi\left([\lambda(\alpha)- \lambda(\Phi(t))]/\tau_1\right) > \Phi\left([\lambda(\alpha)- \lambda(\Phi(t))]/\tau_2\right) = Pr[t < z_2],
$$
and hence $h < 0$, by (\ref{expsd}). Analogously, $Pr[t < z_1]  < Pr[t < z_2]$ if $\Phi(t) - \alpha >0$. It follows that $h < 0$ in each case (and for all $k,\, \theta$), proving that $H_-^s$ holds iff $\tau_1 \le \tau_2$. \done

\section*{Acknowledgments}
This work was funded by the European Union Seventh Framework Programme under grant agreement 290976. The Klaus Tschira Foundation provided infrastructural support at the Heidelberg Institute for Theoretical Studies (HITS). We thank the referees, Tilmann Gneiting, Alexander Jordan, and Wolfgang Polonik for valuable comments and discussions.

\end{document}